\documentclass[11pt, a4paper]{amsart}

\usepackage{graphicx,amsmath,amssymb}
\usepackage{amsthm, mathrsfs}
\usepackage{enumerate}
\usepackage{enumitem}
\usepackage[]{todonotes}
\usepackage[all]{xy,xypic}
\usepackage{url}
\usepackage{xcolor}

\newtheorem{theorem}{Theorem}[section]

\newtheorem{corollary}[theorem]{Corollary}

\theoremstyle{definition}
\newtheorem{definition}[theorem]{Definition}

\newtheorem{lemma}[theorem]{Lemma}

\newtheorem{proposition}[theorem]{Proposition}
\theoremstyle{definition}
\newtheorem{remark}[theorem]{Remark}

\theoremstyle{remark}
\newtheorem{subclaim}{Claim}[theorem]

\newtheorem{pkt}[theorem]{}
\newcommand{\bpk}{\begin{pkt}\rm }  
\newcommand{\epk}{\end{pkt}} 

\newcommand{\satisfies}{\models}
\renewcommand{\bar}{\overline}
\renewcommand{\implies}{\rightarrow}
\renewcommand{\phi}{\varphi}

\newcommand{\ideal}{\triangleleft}
\newcommand{\clsubset}{\subseteq_{cl}}
\newcommand{\opsubset}{\subseteq_{op}}
\newcommand{\pr}{\operatorname{pr}}

\newcommand{\Dom}[1]{\operatorname{Dom}(#1)}

\DeclareMathOperator{\dcl}{\operatorname{dcl}}

\DeclareMathOperator{\id}{\operatorname{id}}

\title[Specialisations of Regular Covers of Zariski Structures]{On the Theory of Specialisations of Regular Covers of Zariski Structures}
\author{U\u gur Efem$^1$}

\author{Boris Zilber$^2$}
\thanks{$^1$Corresponding author. Dyson Institute of Engineering and Technology, ugur.efem@dyson.com}
\thanks{$^2$ University of Oxford, Mathematical Institute.}

\begin{document}
\begin{abstract}
In algebraic geometry specialisations and valuations play 
and important role. In this paper we start investigating 
analogous structures for Zariski structures. Specifically, 
we look into the existence and uniqueness properties 
of extensions of universal specialisations from a base 
Zariski structure to its regular cover. In the process we 
begin to uncover some structural properties of regular 
covers of Zariski structures, and also to uncover the type 
of topological properties necessary for a Zariski structure 
to have a ``good'' theory of specialisations. A subclass 
of Zariski structures is identified with a ``good'' theory 
of specialisations.
\end{abstract}

\maketitle

\section{Introduction}\label{Sec:Intro}
Covers in model theory originated from the study of totally 
and uncountably categorical structures~\cite{AZ86, AZ91, 
Zilber93}. The Ladder Theorem of the second author showed how uncountably 
categorical structures can be built from strongly minimal 
structures by constructing a series of covers. This theorem 
introduced a technique of building new structures from simpler 
ones. Unsurprisingly, some new and non-classical examples 
of Zariski structures are built using ideas of the Ladder Theorem, 
and covers~\cite{SSZ, Zilber08, Zilber10-2}. However, covers 
of Zariski structures are not studied systematically. An initial 
study is started in the first author's PhD thesis~\cite{Efem19}, 
where a certain type of covers where considered. The type of 
covers introduced there includes all the examples of non-classical 
Zariski structures known to us.

On the other hand, specialisations in the setting of Zariski 
structures is an immediate generalisation of specialisations in the sense of A.Weil, 
and valuations from the theory of valued fields. 
They were introduced in~\cite{HZ96} and used to prove the 
Classification Theorem in the same paper. A study of the theory 
of specialisations of Zariski structures was initiated by Onshuus 
and the second author~\cite{OZ}, and motivated by their work, continued 
by the first author in his PhD thesis~\cite{Efem19} under the 
supervision of the second author. The findings of this paper 
are mainly from the last chapter of this thesis. 

Although the theory of specialisations and valuations are very well 
established in algebraic geometry, it is far from being completely 
established generally for Zariski structures. In particular, suppose a Zariski structure 
$\mathcal{C}$
has 
a form of a cover $\mathcal{C}=(C,M, \pr),$ $\pr: C\twoheadrightarrow M,$ and over the Zariski structure $\mathrm{M}$ on $M$ a  specialisation $\pi: \mathrm{M}^*\to \mathrm{M},$ $ \mathrm{M}^*\succ \mathrm{M},$ 
 is given, 
 
 {\em Can one lift $\pi$ to a specialisation over $\mathcal{C}?$}
 
 When $\mathrm{M}$ is just a Zariski structure on an algebraic variety over an algebraically closed field $\mathrm{F},$ a specialisation $\pi$ is given by a valuation over   $\mathrm{F},$
(see more details in ~\cite{Efem}). The theory of the structure $(\mathrm{M}^*, \mathrm{M}, \pi),$ the specialisation theory over $\mathrm{M},$ in this case is bi-interpretable with a well studied theory $\mathrm{ACVF}$ of algebraically closed valued fields, see \cite{HHM1}.  
Note that this theory, by A.Robinson, eliminates quantifiers.

But even in this case the answer to the above question is in general negative.  
This paper identifies a natural condition, 
 the ``Continuous Connections'' 
assumption (CC), which is satisfied when  $\pr: C\twoheadrightarrow M$ is a fibration in the sense of algebraic geometry and which we prove to be sufficient to answer the question in positive.  

Our main result is Theorem~\ref{thm:QE}: 
{\em Assume that the specialisation theory $\mathrm{Th}(M^*,M,\pi_M)$ of the base sort eliminates quantifiers. Then $\mathrm{Th}(\mathcal C)^{\pi},$  the theory of specialisation of the cover structure satisfying (CC),
admits quantifier elimination and  
is complete. }

This is a generalisation of the above mentioned Robinson's theorem.



Now a few words on the assumption (CC). It is not very hard to check, see   that
the example of a non-classical Zariski structure, introduced in~\cite{HZ96}, section 10 
satisfies this assumption.   

The non-classical Zariski structures constructed as covers 
in~\cite{SSZ, Zilber08, Zilber10-2} 
are  objects of non-commutative 
geometry. In contrast to (commutative) algebraic geometry, valuation theory in non-commutative geometry is far from being straightforward. We hope the paper is a contribution towards developing such a theory. 
 An immediate task 
in this direction 
is to establish that all the 
structures satisfy ``Continuous Connections'' assumption for which one requires an efficient enough definability analysis.


\medskip

In Section~\ref{Sec:prelim}, we introduce the preliminary concepts. 
It provides the necessary technical background on Zariski structures, 
their elementary extensions, quotient objects in Zariski structures, 
Zariski groups, and specialisations. The majority of the material 
here is known well are readily available in the model theory literature, 
(~\cite{Efem19, HZ96, Junker, Junker96, Zilber} to name a few). 
However Subsection~\ref{subsec:special} on specialisations, along 
known results, also contains some new technical results and a 
new characterisation for $\aleph_0$-universal specialisations 
for a certain type of structure which is compatible with the regular 
covers we consider in this paper.

Section~\ref{sec:Reg-cvrs} introduces the notion of regular cover 
of a Zariski structure and discusses their structure. We discuss the 
case where there is a single definably almost simple structure group 
acting on the fibres. Under additional but natural assumptions a description 
of closed relations between fibres of the cover is achieved. This leads 
to the ``Continuous Connections'' assumption and is essential for 
the results on specialisations of regular covers discussed in the next 
section. 

Section~\ref{sec:spcl-cvrs} looks into the specialisations of regular 
covers of Zariski structures. We prove that a certain class of regular 
covers which satisfy  the assumptions set in the previous section has 
``good'' theory of specialisations. In particular we prove that any 
maximal extension of a universal specialisation of the base Zariski 
structure to cover is also universal. 

The discussions of Section~\ref{sec:spcl-cvrs} led us to present 
in Section~\ref{sec:theory}, a first order theory for a given regular 
cover together with a specialisation. One can see this theory is a 
generalisation of the theory of an algebraically closed field with a 
specialisation presented in~\cite{Efem}. We prove that this theory 
admits quantifier elimination and is complete. 

As already mentioned, we conclude with a discussion of ``Continuous 
Connections'' assumption in the first example of a non-classical Zariski 
structure.

\section{Preliminaries}\label{Sec:prelim}
We start with introducing the necessary background knowledge 
on Zariski structures, their elementary extensions, quotient objects 
in Zariski structures including Zariski groups and specialisations. 
Most of the material in this section is standard and available in 
the literature, however material on topological sorts and Zariski 
groups contain certain mild generalisations, and slightly different 
perspectives, and in fact, the discussion on specialisations does 
contain some new technical results.

\subsection{Zariski Structures and Elementary Extensions} 
 
\begin{definition}\label{def:zar-str}
A \emph{Noetherian Zariski structure} is a tuple $(M,\{\tau_n:n\in
\mathbb N\},\dim)$ where $\tau_n$ is a Noetherian topology on 
$M^n$ (for each $n$), and $\dim$ is a function which associates 
to every constructible subset of $M^n$ a natural number, which 
also satisfy the following axioms: 
\paragraph{\bf Topological:}
\begin{enumerate}
\item The graph of equality is closed.
\item Any singleton in $M^n$ is closed.
\item Cartesian products of closed sets are closed.
\item The image of a closed set under a permutation of coordinates is closed.
\item For $a\in M^k$ and a closed $S\subseteq M^k+l$, the set $S(a,M^l) :\{m\in M^l: 
(a,m)\in S\}$ is closed.
\item[(SP)]Semi Properness: For a closed irreducible 
$S\subseteq M^n$ and a projection $\pr:M^n\to M^m$, there is a proper 
closed subset $F\subset \bar{\pr(S)}$ such that $\bar{\pr(S)}\setminus 
F \subseteq \pr(S)$. 
\end{enumerate}
\paragraph{\bf Dimension:}
\begin{itemize}
\item[(DP)] Dimension of a Point: $\dim(a) = 0$ for all $a\in M$.
\item[(DU)] Dimension of Unions: $\dim(S_1\cup S_2) = \max(\dim(S_1). 
\dim(S_2))$ for closed $S_1$, and $S_2$.
\item[(SI)] Strong Irreducibility: For any $S\clsubset U\opsubset M^n$ and any closed $S_1
\subsetneq S$, $\dim(S_1)<\dim(S)$.
\item[(AF)] Addition Formula: For any irreducible closed $S\clsubset U
\opsubset M^n$ and a projection $\pr:M^n\to M^m$,
\[\dim(S) = \dim(\pr(S)) + \min_{a\in\pr(S)}(\pr^{-1}(a)\cap S). \]
\item[(FC)] Fiber Condition: Given $S\clsubset U\opsubset M^n$ and a 
projection $\pr:M^n\to M^m$, there is a relatively open $V\opsubset\pr(S)$ 
such that, for any $v\in V$
\[\min_{a\in\pr(S)}(\dim(\pr^{-1}(a)\cap S)) = \dim(\pr^{-1}(v)\cap S).\] 
\end{itemize}

The pair $(M,\{\tau_n: n\in\mathbb N\})$ (without the dimension 
$\dim$) where $\{\tau_n: n\in\mathbb N\}$ satisfies all of the 
topological axioms above is called a \emph{Noetherian topological 
structure}.
\end{definition}

\begin{remark}
Although there is a notion of analytic Zariski structure~\cite[Chap. 6]{Zilber}, 
where Noetherianity is not necessary, this paper is entirely in the context 
of Notherian Zariski structures. Therefore, as there is no danger of ambiguity, 
we will often drop Noetherian, and simply say Zariski structure (or occasionally 
topological structure).
\end{remark}

We can see any Zariski structure $(M,\{\tau_n:n\in\mathbb N\},\dim)$ as 
a first order structure, by introducing a predicate for each closed subset 
in the topology $\tau_n$, for all $n$. Now, let $\mathcal L$ be the first 
order language consisting of the predicates we introduced. Further, we 
assume that $\mathcal L$ contains a constant symbol for each element 
of $M$. Then $(M,\{\tau_n:n\in\mathbb N\})$ becomes an $\mathcal L$ 
structure. So, for all $n\in \mathbb N$, closed subsets of $M^n$ are given 
by positive quantifier free $\mathcal L$-formulas. With a slight abuse of 
the terminology, we will call any element of the set $\cup\{\tau_n:n\in\mathbb N\}$, 
a closed set (of $M$). 

We also immediately see that constructible sets are $\mathcal L$-definable. 
In fact, the converse is also true. Definable sets are constructible. In model 
theoretic terms this is to say $(M,\{\tau_n:n\in\mathbb N\},\dim)$ admits 
quantifier elimination (see~\cite[Theorem 3.2.1]{Zilber}).

In the rest of this paper, as it will not cause any ambiguity, we will write 
$M$ is a Zariski structure instead writing the whole tuple $(M,\{\tau_n:n\in
\mathbb N\}, \dim)$. We will always consider a Zariski structure $M$ as 
a first order structure for some appropriate language as described above. 
Also, for any closed subset $S\in\tau_n$ of $M^n$ we will identify $S$ 
with the positive quantifier free $\mathcal L$-formula which defines it, 
and we will also denote this formula by $S$. Following this convention 
$M\models S(\bar a)$ will mean that $a\in S\subseteq M^n$.

\begin{remark}
A subset $S\subseteq M^n$ is said to be \emph{irreducible} if it cannot 
be written as a union of two proper (relatively) closed subsets. 

As a consequence of Noetherianity, any closed set $S$ can be written as 
a finite union of distinct relatively closed and irreducible subsets uniquely up to 
ordering. They are called \emph{irreducible components} of $S$.
\end{remark}

\begin{remark}\label{rmk:cnstrctbl}
From the definition one can easily observe that a constructible set $Q$ 
can be written as 
\[Q = \bigcup_{i\leq k} S_i\setminus P_i\]
for some $k\in\mathbb N$, and closed sets $S_i,P_i$ such that $P_i\subset S_i$ 
and $S_i$ irreducible. 
Therefore clearly,
\[\bar{Q}=\bigcup_{i\leq k}S_i.\]
\end{remark}

By quantifier elimination, and Remark~\ref{rmk:cnstrctbl}, allows us to extend 
the dimension to definable sets as follows:
\[\dim(Q) = \dim(\bar{Q}) = \max_{i\leq k}\dim(S_i)\]

\begin{definition}
A Zariski structure $M$ is said to be 
\begin{enumerate}
\item \emph{complete} if for any closed 
$S\subseteq M^n$, and any projection $\pr_{i_1,\ldots,i_m}: M^n\to M^m$, 
the set $\pr_{i_1,\ldots,i_m}(S)\subseteq M^m$ is closed. 
\item \emph{pre-smooth} if for any closed irreducible $S_1, S_2\subset M^n$ 
any irreducible component $S$ of $S_1\cap S_2$ satisfies 
\[\dim(S) \geq \dim(S_1) + \dim(S_2) - n\]   
\end{enumerate}
 
\end{definition}

We will work with elementary extensions of Zariski structures. The situation 
with elementary extension is briefly as follows, a more detailed explanation 
can be found in~\cite[Subsection 3.5.3]{Zilber}.

Let $M_0$ be a Zariski structure in an appropriate language $\mathcal L$. 
Let $M\succeq M_0$ be an elementary extension. We will define a topology 
on each $M^n$ as follows: For any closed $S \subseteq M_0^{l+n}$ for 
any $l\in\mathbb N$ declare the subsets of the form $S(\bar a, M^n)$ 
of $M^n$ closed where $\bar a\in M^l$. It is easy to check that this 
collection of subsets is a topology on $M^n$, let us denote it by $\rho_n$. 

\begin{definition}
Let $S\subseteq M_0^{l+n}$ be a closed subset, and $M\succeq M_0$. 
The closed subsets of $M^n$ of the form $S(\bar a, M^n)$ where the \
parameter $\bar a$ is in a subset $A\subseteq M^l$ are called 
\emph{$A$-closed subsets}. 
\end{definition}

We also define a dimension function on the constructible sets in $M$. 
Let $S\subseteq M_0^{l+n}$ be an $M_0$-closed set. Define 
\[\mathcal P(S,k) :=\{\bar a\in \pr(S) : \dim(S(\bar a, M_0^n)) > k\}\]
for the projection $\pr: M_0^{l+n}\to M_0^l$, and where $\dim$ in the 
definition of $\mathcal P(S,k)$ is the dimension function of the structure 
$M_0$. By (AF), dimension of the fibres of $S$ is bounded. So, for every 
$\bar a\in\pr(S)$ there is a maximal $k$ such that $a\in\mathcal P(S,k)$. 
We define a dimension function for the structure $M$ as  
\[\dim(S(\bar a, M^n)) := \max\{k\in\mathbb N: \mathcal P(S,k)\} + 1\]
 
It is immediate from the construction that $(M,\{\rho_n\},\dim)$ satisfies 
the topological axioms. However, ensuring the Noetherianity of the 
topologies $\rho_n$ in general would require a further technical assumption 
that $M_0$ satisfies the following axiom: 
\begin{itemize}
\item[(EU)] Essential Uncountability: Let $S\subseteq M_0^n$ be a 
closed set. If $S$ can be written as a union of countably many closed 
subsets of $M_0^n$, then it can be written as a union of finitely many 
of those subsets.
\end{itemize}
Then, assuming (EU), Noetherianity of the topologies $\rho_n$ 
follows~\cite[Lemma 3.5.24]{Zilber}.

The dimension axioms, (DP), (FC),(DU), (SI) and (FC) are satisfied 
by $(M,\{\rho_n\},\dim)$ Further, if $M_0$ is pre-smooth then the 
elementary extension $M$ is also pre-smooth. Proofs are given 
in detail in~\cite[Subsection 3.5.3]{Zilber}. 

\subsection{Topological Sorts and Many Sorted Zariski Structures}

In general Zariski Structures do not admit elimination of imaginaries 
(see~\cite[p.112]{Marker} for an example,~\cite[sct. 6]{Marker} gives 
more sophisticated examples). However, luckily there are certain 
important sorts one frequently ends up considering which are easy 
to describe and share many important properties of Zariski structures.

Let $D\subseteq M^n$ be a definable subset, and $E\subseteq D\times 
D$ a closed equivalence relation on $D$. Then $T:=D/{E}$ is an 
imaginary sort of $M$. Let $\mathrm{p}:D\to D/{E}$ denote the canonical 
quotient map. When dealing with more than one sort at the same time 
we will add a subscript and write this map as $\mathrm{p}_T$.

Cartesian powers $T^m = (D/{E})^m$ can be identified with $D^m/{E^m}$ 
by declaring $(a_1,\dots,a_m)E^m(b_1,\ldots,b_m)$ if and only if $a_iEb_i$ 
for all $i$. The map $\mathrm{p}:D\to D/{E}$ applied coordinatewise to 
$D^m$ induces a quotient map $D^m\to D^m/E^m$. We will again denote 
this map by $\mathrm{p}$. We equip $T$ and all its Cartesian powers $T^m$ 
with the corresponding quotient topologies via the canonical quotient maps 
from $D$ and $D^m$. 

\begin{definition}
Let $T:= D/{E}$ be an imaginary sort where $D\subseteq M^n$ is definable 
and $E\subseteq D\times D$ is a closed equivalence relation. Then $T$ 
together with the collection of closed subsets of its Cartesian powers is 
called a \emph{topological sort}.
\end{definition} 

\begin{remark}
Originally, in~\cite{Zilber}, topological sorts are defined for $D$ an irreducible 
definable subset of $M^n$. Hence they are irreducible. Here we are extending 
the notion to allow reducible topological sorts, although the only reducible 
topological sorts we will consider in this paper are the finite ones.

\end{remark}

Let $T=D_1/{E_1}$ and $T_2=D/{E_2}$ be two topological sorts in a Zariski 
structure $M$. Then $(D_1\times D_2)/{E_1\times E_2}$ is also a topological 
sort, and it is identified with $T_1\times T_2$.

\begin{proposition}
A topological sort $T$ satisfies all of the topological axioms (given in 
Definition~\ref{def:zar-str}). I.e. $T$ is a topological structure with 
the induced quotient topology.
\end{proposition}

\begin{proposition}\label{prop:top-Meq}
Let $T=D/{E}$ be a topological sort. Let $E^\prime$ be a closed equivalence 
relation on $T$. Then $T/{E^\prime}$ is a topological sort, and it can be 
represented as $D/{E^{\prime\prime}}$ for some closed equivalence relation 
$E^{\prime\prime}\subseteq D^2$.
\end{proposition}
\begin{proof}
Let $T=D/{E}$ be a topological sort and $E^\prime\subseteq T^2 = D/{E}
\times D/{E}$ be a closed equivalence relation. Define $E^{\prime\prime}$ 
as $\mathrm{p}^{-1}(E^\prime)$. Clearly $E^{\prime\prime}\subseteq D\times D$ 
is closed. It is easy to see $E^{\prime\prime}$ is an equivalence relation. 
Let $d\in D$, and by $[d]$ denote the $E$-equivalence class of $d$. Since 
$E^\prime$ is an equivalence relation on $T$, we see $[d]E^\prime [d]$. 
Hence $(d,d)\in\mathrm{p}_T^{-1}(d)\times\mathrm{p}_T^{-1}(d)\subseteq 
E^{\prime\prime}$. It follows, via similar arguments, that $E^{\prime\prime}$ 
is symmetric an d transitive. 
\end{proof}

\begin{definition}
Let $T = D/{E}$ be a topological sort, and $F\subseteq T^n$ be a closed and 
irreducible subset with $\mathrm{p}^{-1}(F) = S$. Then we define dimension 
of $F$ by 
\[\dim(T) := \dim(S) -\min\{\dim(\mathrm{p}^{-1}(f)) : f\in F\}\]

For any arbitrary closed subset, we define the dimension to be the maximum 
of dimensions of its irreducible components.
\end{definition}


\begin{lemma}
Let $T$ be a topological sort with the dimension $\dim$ defined as above. 
Then $(T,\dim)$ satisfies the dimension axioms (DP), (DU), and (SI).
\end{lemma}
\begin{proof} Immediate from definitions.
\end{proof}
\begin{definition}
Let $M$ be a Zariski structure and $T_1, T_2$ be two topological sorts of 
$M$. A function $f:T_1\to T_2$ is called a morphism if $f\times \id: T_1
\times M^n \to T_2\times M^n$ is continuous for all $n$.
\end{definition}

\begin{lemma}\label{lem:morph-prop-1}
\begin{enumerate}[label=(\roman*)]
\item Graphs of morphisms are closed.
\item For a topological sort $T=D/{E}$, the quotient map $\mathrm{p}_T:
D\to T$ is a morphism.
\item If $f:T_1\to T_2$, $g:T_2\to T_3$ and $h:T_3\to T_4$ are morphisms, 
then $g\circ f:T_1\to T_2$ and $f\times h:T_1\times T_3 \to T_2\times T_2$ 
are morphisms.
\end{enumerate}
\end{lemma}
\begin{proof} Immediate from definitions.
\end{proof}
\begin{lemma}\label{lem:morph-prop-2}
Let $f:T_1\to T_2$ be a morphism. 
\begin{enumerate}[label=(\roman*)]
\item Let $E_f$ be the equivalence relation given 
by the pre-image of $=$ under $f$. Then $f$ factorises as
\[\xymatrix{ T_1 \ar@{->>}[r]^{p} &T_1/{E_f}\ar@{->}[r]^{\tilde{f}} & Im(f) \ar@{^{(}->}[r]^{i}& T_2}\]
where $\tilde{f}$ is a bijective morphism and $i$ is the inclusion map.
\item f $V\subseteq T_1$ is definable and irreducible, then $f_{|_V}:V\to T_2$ 
is a morphism. 
\item If $E$ is a closed equivalence relation on $T_1$ and $f$ is constant 
on its equivalence classes, then $f$ induces a morphism $f/{E}:T_1/{E}
\to T_2$.
\end{enumerate}
\end{lemma}
\begin{proof} Immediate from definitions.
\end{proof}
\begin{corollary}
Let $T_1$ and $T_2$ be topological sorts, and $f: T_1\to T_2$ a surjective 
morphism. Let $E(x_1,x_2)$ be the equivalence relation $f(x_1)=f(x_2)$ on 
$T_1$. Then $T_1/E$ is a topological sort and $f$ induces a homeomorphism 
between $T_1/E$ and $T_2$.
\end{corollary}

\begin{lemma}
Let $T_1$ and $T_2$ be topological sorts, then $f:T_1\to T_2$ is a 
morphism if and only if for any $n$ and any positive quantifier free 
formula $\psi(x,\overline{z})$ defining a closed subset of $T_1\times 
M^n$, the formula $\exists y\; f(x) = y \wedge \psi(y,\overline{z})$ 
is equivalent to a positive quantifier free formula. We will 
denote this formula by $\psi(f(x),\overline{z})$. 
\end{lemma}
\begin{proof}
Assume $f$ is a morphism. The formula $\exists y\; f(x) = y \wedge 
\psi(y,\overline{z})$ defines the pre-image of a closed set under the 
map $f\times id: T_1\times M^n \to T_2\times M^n$. Since $f$ is a 
morphism, this map is continuous. Hence the pre-image is closed. 
Which means $\exists y\; f(x) = y \wedge \psi(y,\overline{z})$ is 
equivalent to a positive quantifier free formula. 
Conversely, the assumption immediately implies that $f\times \id: 
T_1\times M^n\to T_2\times M^n$ is continuous for all $n$. Hence 
$f$ is a morphism.
\end{proof}

\begin{definition}
Let $H$ be a Zariski structure. A multi-sorted Zariski structure $M$ with 
the home sort $H$ is a multi-sorted structure in a multi-sorted language 
$\mathcal L$ with sorts $(M_i)_{i\in I}$ such that: 
\begin{enumerate}[label=(\roman*)]
\item Each $M_i$ is a topological sort in $H$; and there is an $i$ such that 
$M_i = H$.
\item Each $M_i$ is an $\mathcal L_i$-Zariski structure where $\mathcal L_i$ 
is the natural language for the sort $M_i$ and $\mathcal L_i\subset\mathcal L$.
\item If $M_i$ and $M_j$ are sorts in $M$ then their product $M_i\times M_j$ 
is also a sort in $M$.
\end{enumerate}
If all sorts $M_i$ of $M$ are pre-smooth, then $M$ is called a \emph{multi-sorted 
Zariski geometry} 
\end{definition}

It is important to remark that, a multi-sorted Zariski structure $M$ is not 
necessarily a Zariski structure with the given language $\mathcal L$. Having 
a home sort $H$ is important, as all other sorts are topological sorts in $H$, 
dimension on each sort $M_i$ is induced by the dimension of $H$ as explained. 
The structure $M$ is ``closed under products of sorts'' in the sense that if 
$M_i$ and $M_j$ are sorts in $M$, then so is $M_i\times M_j$. 

\begin{proposition}
Let $M$ be a multi-sorted Zariski structure with a home sort $H$. Then any 
sort $M_i$ of $M$ is stably embedded.
\end{proposition}
\begin{proof}
Since $H$ is stable, $H^{eq}$ is also stable. Therefore any sort interpretable 
in $H$ is stably embedded. 
\end{proof}

\subsection{Zariski Groups}

\begin{definition}\label{def:z-topl-grp}
Let $C$ be a Zariski structure, and $G$ be a group that is a topological sort 
in $C$ such that multiplication $m:G\times G\to G$ and inversion ${}^{-1}:G
\to G$ are morphisms. Such a group $G$ will be called a \emph{Zariski group 
(in $C$)}. 
\end{definition} 

\begin{remark}\label{rmk:non-conn-finite-grps}
As we are allowing reducible definable sets in their construction, topological 
sorts in general are not necessarily irreducible. For Zariski groups this in particular 
will allow us to consider non-connected (stable) groups as Zariski groups. 
Although, in this paper the only non-connected Zariski groups we will consider 
are the finite ones. Here we only present some essential properties of Zariski 
groups we will require in this paper. For a more detailed study of these groups 
we refer the reader to~\cite{Junker, Junker96}.
\end{remark}

\begin{definition}\label{def:cont-prop-act}
Let $C$ be a Zariski structure and $G$ be Zariski topological group in $C$.
Let $A\subseteq C$ be a constructible set. 
\begin{enumerate}[label=(\roman*)]
\item
We say that \emph{$G$ acts morphically on $A$} (or the action of $G$ is morphic)
if the action 
\begin{eqnarray*}
\Theta: G\times A &\to& A\\
(g,a) &\mapsto & g\cdot a
\end{eqnarray*} 
is a morphism. Often we will denote the action by $\cdot$, and write $g\cdot a$.
\item We say that the action is \emph{proper}
if
\[E_G(x,y) \mbox{ defined by } \exists g\in G\,(y = g\cdot x) \]
is a closed equivalence relation in $C\times C$.
\item We will say that the action of $G$ is \emph{free} 
if it is proper, and 
the action $G\times A\to A$ is invertible in the sense that  
there is a morphism
\begin{eqnarray*}
&E_G &\to  G\\
&(a,g\cdot a) &\mapsto  g
\end{eqnarray*}
\end{enumerate}
\end{definition}

In this paper, whenever we talk about the action of a Zariski group, 
we will always assume that the action is morphic.

\begin{proposition}[Proposition 5.6 in~\cite{Junker96}]\label{fct:clsr-subgr-constr-mono}
Let $G$ be a Zariski group, and $H< G$ be a subgroup (not necessarily definable!). 
Then the closure $\overline H$ is also a subgroup. If $H$ is normal, then so is 
$\overline H$. Moreover, a definable submonoid of $G$ is a closed subgroup.
\end{proposition}
\begin{proof}
Let $H< G$ be a subgroup. Then $H\times H\subseteq m^{-1}(\overline H)$ 
where $m:G\times G\to G$ is the group operation. Since $m$ is continuous, 
$m^{-1}(\bar H)$ is closed. Then it follows 
\[\overline H\times \overline H = \overline{H\times H} \subseteq m^{-1}(\overline H)\]
Hence 
\[m(\overline H\times \overline H) = \overline H \cdot \overline H\subseteq \overline H\]
 
A similar argument with inversion instead of multiplication shows $\overline H$ 
is a subgroup.

Let $g\in N_G(H)$. Conjugation with $g$ is a homeomorphism, so $\overline H^g = 
\overline{H^g} = \overline H$. Then $N_G(H)\subseteq N_G(\overline H)$.

By stability, any definable submonoid is a subgroup. All cosets of the subgroup 
$H$ in $\overline H$ are homeomorphic to $H$. Hence all are dense in $\overline H$. 
But two disjoint constructible sets cannot be dense in the union. Hence $H=\overline{H}$.
\end{proof}

\begin{lemma}\label{lem:z-topl-grp-quot}
Suppose $G$ is a Zariski group acting freely on $D$ and $H \ideal G$ 
is a definable normal subgroup. Then $G/{H}$ is a Zariski group, $D/{H}$ 
is a topological sort, and $G/{H}$ acts freely on $D/{H}$. 
\end{lemma}
\begin{proof}
We can assume $H$ is closed (see Fact~\ref{fct:clsr-subgr-constr-mono}). 
First let us show that the quotient $D/{H}$ is a topological sort. Since the action 
is free there is a morphism $E_G\to G$ defined by $(a, g\cdot a)\mapsto g$. 
Let $E_H$ be the pre-image of $H$ under this map. Clearly, $E_H$ is a definable 
equivalence relation on $D$ given by 
\[x E_H y :\Leftrightarrow \exists h\in H\, (h\cdot x = y)\]
 
Since $H$ is closed its inverse image $E_H$ under this morphism is closed in 
$E_G$. Hence $E_H$ is a closed equivalence relation on $D$. Therefore $D/{H}$ 
is a topological sort. 

Next, we sill how that $G/{H}$ is a topological group. The closed normal subgroup 
$H$ defines the equivalence relation
\[a\sim b :\Leftrightarrow ab^{-1}\in H\]
for all $a,b\in G$. Observe that $\sim$ is the pre-image of $H$ under 
the morphism
\begin{eqnarray*}
G\times G &\to& G\\
(a,b) &\mapsto & ab^{-1}
\end{eqnarray*}  
Since $H$ is closed, $\sim$ is closed. Hence $G/{\sim} = G/{H}$ is 
a topological sort in $G$. To see that multiplication on $G/{H}$ is a 
morphism observe that $\mathrm{p}_H\circ m: G\times G\to G \to 
G/{H}$ is a morphism. Moreover $H\times H$ is a closed equivalence 
relation on $G\times G$, such that the morphism $\mathrm{p}_H\circ 
m$ is constant on its classes. By Lemma~\ref{lem:morph-prop-2}, 
multiplication is a morphism. A similar argument will show that inversion is also 
a morphism. 

By Proposition~\ref{prop:top-Meq}, $G/{H}$ is also a topological sort in $C$, 
and the group operations are again morphisms.

Next, let us consider the action of $G/{H}$ on $D/{H}$. The action $\Theta: G 
\times D\to D$ and the quotient map $\mathrm{p}_H:D\to D/{H}$ are morphisms. 
Therefore their composition 
\[\mathrm{p}_H\circ\Theta:G\times D\to D\to D/{H}\]
is a morphism. Moreover $\mathrm{p}_H\circ\Theta$ is constant on equivalence 
classes of the closed equivalence relation defined by $H\times H$ on $G\times 
D$. Therefore 
\[\mathrm{p}_H\circ\Theta/{H\times H}:G\times D/{(H\times H)}\to D/{H}\]
is a morphism. Since $G\times D/({H\times H)}\simeq G/{H}\times D/{H}$, 
the action of $G/{H}$ on $D/{H}$ is a morphism. 
\end{proof}

\begin{lemma}\label{lem:quot-fin-z-top-grp}
Let $H$ be a finite Zariski group acting freely on $D$ and $T = D/{H}$ 
be a topological sort, let $\mathrm{p}_T:D\to T$ be the canonical quotient 
map. Let $Q\subseteq D$ be a closed subset. Then $\mathrm{p}_T(Q)
\subseteq T$ is also closed.
\end{lemma}
\begin{proof}
 Note that
\[\mathrm{p}_T^{-1}(\mathrm{p}_T(Q))=H\cdot Q= \bigcup_{h\in H} h\cdot Q\]
Since subsets $h\cdot Q$ are closed, the statement follows.
\end{proof}

\subsection{Specialisations}\label{subsec:special}
Specialisations are the main objects we study in this paper. Here 
we describe specialisations, briefly explain why they are important 
to study, and layout important properties which we will use later 
in the paper, where we study specialisations on regular covers in 
detail.

\begin{definition}\label{def:specialization}
Let $M_0$ be a Zariski structure and $M\succeq M_0$. A partial 
function $\pi:M\to M_0$ such that 
\begin{enumerate}[label=(\roman*)]
\item $\pi(m) = m$ for all $m\in M_0$;
\item for every formula $S(\bar{x})$ over $\emptyset$, defining an 
$M_0$-closed set and for every $\bar{a}\in M^n\cap (\Dom{\pi})^n$ 
\[M\satisfies S(\bar{a}) \mbox{ implies } M_0\satisfies S(\pi \bar{a}) \]
\end{enumerate} 
is said to be a \emph{specialisation}. 
\end{definition}

\begin{remark}
Although we defined specialisation for Zariski structures, the definition 
only involves the topological structure $\{\tau_n:n\in\mathbb N\}$ 
(note that the definition above does not involve the dimension function 
$\dim$). So, specialisations can be defined for topological structures 
$(M_0,\{\tau_n: n\in\mathbb N\})$ (observe that an elementary extension 
$M\succeq M_0$ will again be a topological structure).

At various points in the rest of the papers we will consider specialisations 
on topological sorts. Although topological sorts are not necessarily Zariski 
structures, they are topological structures. As explained here topological 
structures are enough to consider specialisations. So we will be safe in 
doing this.
\end{remark}

\begin{definition}
A specialisation is said to be \emph{$\kappa$-universal} if, given any 
$M^\prime\succeq M\succeq M_0$, any $A\subseteq M^\prime$ 
with $|A|<\kappa$ and a specialisation $\pi_A:M\cup A\to M_0$ 
extending $\pi$, there is an elementary embedding $\sigma:A\to 
M$ over $M\cap A$ such that $\pi_A|A = \pi\circ\sigma$.
\end{definition}

\begin{definition}\label{def:max-special}
Let $\pi:M\to M_0$ be  a specialisation. We say that $\Dom{\pi}$ is 
\emph{maximal} if there are no specialisations $\pi^\prime: M\to M_0$ 
extending $\pi$ non-trivially.
\end{definition}

\begin{proposition}\label{prop:max-vs-univ-spcl}
An $\aleph_0$-universal specialisation is maximal.
\end{proposition}
\begin{proof}
Let $M_0\preceq M$ be a pair of Zariski structures and $\pi: M\to 
M_0$ be an $\aleph_0$-universal specialisation. Assume $\pi$ is not 
maximal. Then there is an $m\in M\setminus\Dom{\pi}$ such that there 
is a specialisation $\pi_{\{m\}}:\{m\}\cup \Dom{\pi}\to M_0$ extending 
$\pi$. Since $\pi$ is universal, there is an embedding $\sigma:\{m\}
\to M$ over $\{m\}\cap M$ with ${\pi_{\{m\}}}_{|_{\{m\}}} = \pi\circ
\sigma$. Since $m\in M$, we have $\sigma(m) = m$. Which implies 
$\pi(m) = \pi_{\{m\}}(m)$. In particular it means $\pi$ is already defined 
on $m$. 
\end{proof}


Next we describe the relation between specialisations and the topology 
on a Zariski structure. Most importantly that all the topological data 
of a Zariski structure can be recovered from $\aleph_0$-universal 
specialisations of the structure.

\begin{definition}
Let $\pi:M\to M_0$ be a specialisation, a definable relation $S\subseteq M_0^n$ 
is said to be $\pi$-closed whenever $\pi({}^*S)\subseteq S$ where ${}^*S$ 
is the interpretation of $S$ in $M$. 
\end{definition}

The family of $\pi$-closed sets satisfies the topological axioms (Exercise 2.2.9 
in~\cite{Zilber}). Further, we can characterise closed sets of a Zariski structure 
in terms of $\pi$-closed relations. For a Zariski structure $\mathcal C_0$, 
if a definable relation $T$ is $\pi$-closed for every specialisation $\pi:
\mathcal C\to \mathcal C_0$ of $\mathcal C_0$, then $T$ is positive 
quantifier free. In fact this is an instance of a more general result of van 
den Dries~\cite{vdD}


Universal specialisations provide us a stronger result characterising 
positive quantifier free formulas: If $\pi:\mathcal C\to \mathcal C_0$ 
is an $\aleph_0$-universal specialisation then any definable relation 
$S\subset M_0^n$ is closed if and only if it is $\pi$-closed (Proposition 
2.2.24 in \cite{Zilber}). 

\begin{proposition}\label{prop:spcl-rsrtc-def}
Let $\mathcal C_0\preceq \mathcal C$ be a Zariski structure and its 
extension. Let $\pi:\mathcal C\to \mathcal C_0$ be a specialisation. 
Let $R$ be a $\emptyset$-definable set. Define $\pi_{|_R}:R(\mathcal C) 
\to R(\mathcal C_0)$ as $\pi_{|_R}(r) = \pi(r)$ whenever $\pi(r)\in 
R(\mathcal C))$. 
Then $\pi_{|_R}:R(\mathcal C)\to R(\mathcal C_0)$ is a specialisation.  
\end{proposition}
\begin{proof}
By construction, $\pi_{|_R}$ is only defined on the points of $R(\mathcal C)$ 
whose images under $\pi$ are in $R(\mathcal C_0)$). Let $S$ be a 
closed subset of $R(\mathcal C)^n$ and $\bar a\in\Dom{\pi_{|_R}}$ 
with $R(\mathcal C)\models S(\bar a)$. Let $\overline{S}$ denote the 
closure of $S$ in $\mathcal C^n$. Since $\pi$ is a specialisation, 
$\mathcal C_0\models S(\pi(\bar a))$. Since $\bar a\in R(\mathcal C)^n$, 
we have $\pi(\bar a) = \pi_{|_R}(\bar a)$. Hence $R(\mathcal C_0)\models 
S(\pi_{|_R}(\bar a))$ as required.
\end{proof}

For a specialisation $\pi:\mathcal C\to \mathcal C_0$ and a definable set 
$R$ as in Proposition~\ref{prop:spcl-rsrtc-def}, we call the specialisation 
$\pi_{|_R}$ the \emph{restriction of $\pi$ to $R$}. For the sake of simplicity 
we will omit the subscript and write $\pi:R(\mathcal C)\to R(\mathcal C_0)$ 
when referring the restriction.

In the remaining of this section we discuss tools we will frequently in the 
study of specialisations on regular covers. 

\begin{lemma}\label{lem:ext-spcl-via-morph}
Let $D$ and $R$ be topological sorts in a Zariski structure $\mathcal C$. 
Let $f:D\to R$ be a morphism, and let $\pi_D: D(\mathcal C^\prime) \to 
D(\mathcal C)$ be a specialisation with $\Dom{\pi_D}\subseteq D(\mathcal C^\prime)$. 
Then $\pi_D$ induces a unique extension to $\pi_R:R(\mathcal C^\prime) 
\to R(\mathcal C)$, with $\Dom{\pi_R} = f(\Dom{\pi_D})$, given by 
$\pi_R(f(x)) := f(\pi_D(x))$. Moreover, $\pi_R\times \pi_D: R(\mathcal C^\prime)
\times D(\mathcal C^\prime)\to R(\mathcal C)\times D(\mathcal C)$ 
is also a specialisation.
\end{lemma}
\begin{proof}
Let $\psi(y)$ define a closed subset of $R(\mathcal C^\prime)$ (over $\emptyset$). 
We need to show that for any $x_0\subset \Dom {\pi_D}$, and $y_0=f(x_0)$,
\[\vDash \psi(y_0) \Rightarrow\;   \vDash \psi(y_0^{\pi_R})\]
where $y_0^{\pi_R}= f(x_0^{\pi_D})$.
Equivalently,  
\[\vDash \psi(f(x_0)) \Rightarrow\;    \vDash \psi(f(x_0^{\pi_D}))\]
which follows from the fact that $\psi(f(x))\equiv \exists y\, f(x)=y\wedge \psi(y)$ 
is the pre image of $\psi(y)$ under $f$, hence closed in $D(\mathcal C^\prime)$.

For the moreover part observer that $f\times\mathrm{id}: D\times D\to R\times D$ 
is a morphism. Then one can repeat a similar argument with $f\times \mathrm{id}$ 
in place of $f$.
\end{proof}

\begin{corollary}\label{cor:ext-spc-top-srts}
Let $\pi:M\to M_0$ be a specialisation. Let $T$ be a topological sort 
in $M_0$. Then there is 
a unique specialisation $\pi_T : T(M)\to T(M_0)$ commuting with  $\mathrm{p_T}.$ 
\end{corollary}
\begin{proof}
In Lemma~\ref{lem:ext-spcl-via-morph} take the topological sort $D$ 
to be $M$, so the realisation of $D$ in $M_0$ will be $M_0$. Take $R = T$,
and take the map $f$ to be the canonical quotient map $\mathrm{p_T}$. 
\end{proof}

Corollary~\ref{cor:ext-spc-top-srts} is an important instance of 
Lemma~\ref{lem:ext-spcl-via-morph}, which allows us to extend 
a specialisation to topological sots. Whenever we are given a specialisation 
$\pi:\mathcal C\to \mathcal C_0$, and a topological sort $T$ we 
will assume that $\pi$ is extended to $\pi_{T}:T(\mathcal C)\to 
T(\mathcal C_0)$ via the canonical quotient maps. We will often drop 
the subscript and write $\pi:T(\mathcal C)\to T(\mathcal C_0)$. When  
we are considering the specialisation as extended to certain topological 
sorts will also write 
\[\pi:\mathcal C\cup T_1(\mathcal C)
\ldots \cup T_n(\mathcal C)\to \mathcal C\cup T_1(\mathcal C)\ldots 
\cup T_n(\mathcal C)\] 
where $T_i$ are topological sorts.

\begin{remark}
Another useful instance of Lemma~\ref{lem:ext-spcl-via-morph} is when $R$ 
and $D$ are definable sets in $\mathcal C$.
\end{remark}

It is of course possible that one can interpret the same topological 
sort via different definable sets, and different quotient maps. In the 
next proposition we will show that if one extends a specialisations 
to the same topological sort via two different morphisms, the extensions 
are compatible. 
Consider a Zariski structure and its extension $\mathcal C_0\prec 
\mathcal C$, and a specialisation $\pi: \mathcal C \to \mathcal C_0$. 
Let $T$ be a topological sort, $f: A \to T$ and $g: B \to T$ be two 
morphisms where $A$ and $B$ are definable sets in $\mathcal C$. 
Also consider the restrictions $\pi:A(\mathcal C)\to A(\mathcal C_0)$ 
and $\pi:B(\mathcal C)\to B(\mathcal C_0)$ of $\pi$ to $A$ and $B$. 
Let $\pi_f:\mathcal C\cup T(\mathcal C) \to \mathcal C_0\cup 
T(\mathcal C_0)$ and $\pi_g:\mathcal C\cup T(\mathcal C) \to 
\mathcal C_0\cup T(\mathcal C_0)$ be the specialisations extending 
$\pi$ to $T$ via $f$ and $g$ respectively.

\begin{proposition}\label{prop:cmnext-spcl-to-srts}
On the intersection $\Dom{\pi_f}\cap\Dom{\pi_g}$ of their domains 
$\pi_f = \pi_g$. Moreover $\pi_f\cup \pi_g$ is a specialisation extending 
both. 
\end{proposition}
\begin{proof}
Let $t\in \Dom{\pi_f}\cap\Dom{\pi_g}$. By construction (see 
Lemma~\ref{lem:ext-spcl-via-morph}), there are $d\in f^{-1}(t)\cap
\Dom{\pi}$ and $e\in g^{-1}(t)\cap \Dom{\pi}$. Remark that $f(x) = 
g(y)$ defines a closed subset of $A\times B$. Therefore, since $f(e) = g(d)$ 
and since $\pi$ is a specialisation, $f(\pi(e)) = g(\pi(d))$. That is $\pi_f(t) 
= \pi_g(t)$.

Next, we will show that $\pi^0 := \pi_f\cup \pi_g$ is a specialisation 
which is a common  extension. Let $S$ be a closed subset of $\mathcal C^n
\times T^k\times T^l$, and $\bar z\subset\Dom{\pi^0}$ be such that $\models 
S(\bar z)$. We may assume that $\bar z = z_1z_2z_3$ is partitioned 
such that $z_1\in \mathcal C^n\cap(\Dom{\pi})^n$, $z_2\in T^k\cap(\Dom{\pi_f})^k$ 
and $z_3\in T^l\cap(\Dom{\pi_g})^l$.

Define 
\[C:=\{(x_1,\ldots,x_n,y_1\ldots,y_k,t_1,\ldots,t_l)\in\mathcal C^{n+k+l} 
: S(\bar x, f(\bar y), \mathrm{p}(\bar t))\}\]
Clearly $C$ is a closed subset. Moreover, Since $\models S(\bar z)$, 
and $\bar z \subset\Dom{\pi^0)}$, there is a tuple $x_1,\ldots,x_n,d_1,
\ldots,d_k,e_1,\ldots,e_l\subset\Dom{\pi}$ such that $\models C(\bar x,
\bar d,\bar e)$ and $f(\bar d) = z_2$ and $g(\bar e) = z_3$. 
Since $C$ is closed $\models C(\pi(\bar x),\pi(\bar d),\pi(\bar e))$. Then 
by definition $\models S(\pi(\bar x), f(\pi(\bar d)), g(\pi(\bar e))$. 
Also by construction $f(\pi(\bar d)) = \pi_f(z_2)$ and $g(\pi(\bar e)) 
= \pi_g(z_3)$. Hence $\models S(\pi(z_1),\pi^0(z_2),\pi^0(z_3))$.
\end{proof}

\begin{lemma}\label{lem:ext-spcl-to-sort-fin-grp}
Let $\mathcal C_0$ be a Zariski structure, and $\mathcal C\succeq 
\mathcal C_0$ be an extension. Let $T = D/H$ be an orbifold where $H < G$ 
is a finite $\emptyset$-definable subgroup acting freely on a constructible 
set $D\subset\mathcal C$. Let $\pi:\mathcal C\cup T(\mathcal C)\to 
\mathcal C_0\cup T(\mathcal C_0)$ be a maximal specialisation, and let 
$t\in T(\mathcal C)\cap \Dom{\pi}$. Then $\mathrm{p}^{-1}(t) \subset\Dom{\pi}$.
\end{lemma}
\begin{proof}
Consider an $a\in \mathrm{p}^{-1}(t)$ and a positive quantifier free formula 
$Q(y,z)$, such that for some $c\subset C$, the formula $Q(y,c)$ is the locus 
of $a$ over $\Dom{\pi}$. Denote 
\[Q(x,y,z):\equiv Q(y,z) \ \&\ \mathrm{p}_T(y)=x\]

By Lemma~\ref{lem:quot-fin-z-top-grp} we see that the formula $\exists y\, 
Q(x,y,z)$ defines a closed subset which by construction contains $(t,c)$ and 
so does contain $(\pi(t),\pi(c))$. The latter means that 
\[\mathcal C_0\vDash \exists y (Q(y,\pi(c))\ \&\ \mathrm{p}_T(y)=t)\]
Let $a_0$ satisfy the formula $Q(y,\pi(c))\ \&\ \mathrm{p}_T(y)=t$.

Now it is clear that setting $\pi(a):=a_0$ we will have extension of the specialisation 
$\pi$ to $a$. Since $\pi$ is maximal, $a$ must be in $\Dom\pi$.

Since $H$ is a substructure of the prime model, $H\subset \Dom\pi$. Recall 
that $\mathrm{p}^{-1}(t) = H\cdot a$.
By Lemma~\ref{lem:ext-spcl-via-morph}, $H\cdot a\subset \Dom\pi$.
\end{proof}

\begin{lemma}\label{lem:prime-&-mnml}
Let $\mathcal C$ be a Zariski structure which is prime and minimal over 
an $\emptyset$-definable subset $M$. Then $\mathcal C$ is atomic over 
$M\cup A$ for any $A\subseteq C.$
\end{lemma}
\begin{proof}
The theory is $Th(\mathcal C)$ is $\omega$-stable. Then there exists 
$\mathcal C^\flat\succeq\mathcal C$ which contains $M\cup A$ and is atomic 
over $M\cup A$.

Since $\mathcal C$ is prime and minimal over $M$, we get $\mathcal C^\flat = 
\mathcal C$. Hence $\mathcal C$ is atomic over $M\cup A$.
\end{proof}

Following Lemma~\ref{lem:prime-&-mnml}, a useful characterisation of 
$\aleph_0$-universal specialisations can be given. 
Note that the second condition of the Theorem~\ref{thm:charac-univ-spc} 
below implies ``$\pi$ is an $\aleph_0$-universal specialisation''. 
 
\begin{theorem}\label{thm:charac-univ-spc}
Let $\mathcal C_0\prec \mathcal C$ be a Zariski structure and its $\aleph_0$-saturated 
extension, $\pi: \mathcal C\to \mathcal C_0$ a specialisation. Suppose 
also that there is a $\emptyset$-definable subset $M$ such that for every 
$A\subset C$, $\mathcal C$ is atomic over $M\cup A$. 

Then the following are equivalent:

\begin{enumerate}[label=(\roman*)]
\item For any finite $c\subset M\cup \Dom{\pi}$, any finite tuple $b^\prime
\in\mathcal C^\prime\succeq \mathcal C$, and for any specialisation 
\[\pi^\prime: Cb^\prime\to C_0\]
extending $\pi$, there is $b\subset \Dom{\pi}$ such that $b\equiv_c b^\prime$ 
and $\pi(b)=\pi^\prime(b^\prime)$.
\item For any finite $a\subset \mathcal C$, any finite tuple $b^\prime
\in\mathcal C^\prime\succeq \mathcal C$, for any specialisation 
\[\pi^\prime: Cb^\prime\to C_0\]
extending $\pi$, there is $b\subset \Dom{\pi}$ such that $b\equiv_ab^\prime$ 
and $\pi(b)=\pi^\prime(b^\prime)$.
\end{enumerate}
\end{theorem}
\begin{proof}
We only need to prove (i) implies (ii). Suppose for a contradiction that (ii) 
fails for some $\pi^\prime$, $a$ and $b^\prime$. Then, $a$ cannot be 
a subset of $\Dom\pi$, as we are assuming (i). 


Let $p = \mathrm{tp}(a/M\cup\Dom \pi)$. By assumption, $\mathcal C$ 
is atomic over $M\cup\Dom{\pi}$. Therefore $p$ is principal; so it is 
equivalent to a formula $P(z)$ over some $c\subset M\cup \Dom{\pi}$.

Let $q(z,y) = \mathrm{tp}(a,b^\prime/\emptyset)$ and 
\[t(y)= \bigwedge_{Q\in q}\exists z\ P(z) \ \&\ Q(z,y)\]

By construction, $t$ is a type over $c$. Clearly, $b^\prime$ realises $t$. 
By (i) there is $b\in \Dom \pi$ realising $t$, and $\pi(b)=\pi^\prime(b^\prime)$. 
Then $\{ P(z)\} \cup q(z,b)$ is consistent, and by saturation must have a 
realisation in $\mathcal C$. Since $P$ is complete (is an atom) over 
$M\cup\Dom \pi$, we have $P(z)\vdash q(z,b)$. It follows $\vDash q(a,b)$, 
and  $b\equiv_ab^\prime$. A contradiction to our assumptions.
\end{proof}

\subsubsection{Specialisations in Many Sorted Zariski Structures}
Here we briefly discuss the notion of specialisation in many sorted 
Zariski structures. Everything we discussed above will still be valid, 
however one needs to introduce specialisations for many sorted 
Zariski structures rigorously for the sake of completeness. The notion 
is in fact subtly hinted in Lemma~\ref{lem:ext-spcl-via-morph}. 

The specialisation we consider in Section~\ref{sec:spcl-cvrs} are 
actually in the many sorted setting, although it is rather implicit 
and the specialisations in question could easily be considered as 
an extension of a specialisation from the home sort to a topological 
sort as in Lemma~\ref{lem:ext-spcl-via-morph}, we still introduce 
the notion to give a more complete picture.

\begin{definition}
Let $M$ be a multi-sorted Zariski structure with sorts $S$. Let $N$ 
be an elementary extension of $M$. A map $\pi = (\pi_{s_1},
\ldots,\pi_{s_n}):N_{s_1}\times\ldots\times N_{s_n}\to M_{s_1}\times
\ldots\times M_{s_n}$ is called a specialisation if each $\pi_{s_i}:N_{s_i}
\to M_{s_i}$ is a specialisation.
 
It is said to be \emph{$\kappa$-universal}
if, given any $N^\prime\succ N\succ M$, any $A_{s_1}\subseteq N^\prime_{s_1},\ldots,
A_{s_n}\subseteq N^\prime_{s_n}$ with $|A_{s_i}|<\kappa$ for each $i$ and a specialisation 
$\pi_A = (\pi_{A_{s_1}},\ldots,\pi_{A_{s_n}}): N_{s_1}\cup A_{s_1}\times\ldots\times N_{s_n}\cup 
A_{s_n}\to M_{s_1}\times\ldots\times M_{s_n}$ extending $\pi= (\pi_{s_1},\ldots,\pi_{s_n})$, 
there is an embedding $\sigma =(\sigma_{s_1},\ldots,\sigma_{s_n}):A_{s_1}\times\ldots\times 
A_{s_n}\to M$ over $(N_{s_1}\times\ldots\times N_{s_n})\cap (A_{s_1}\times\ldots\times A_{s_n})$ 
such that $\pi_{A_{s_i}}|A_{s_i} = \pi_{s_i}\circ\sigma_{s_i}$ for each $i$.  
\end{definition}  

\section{Regular Covers of Zariski Structures}\label{sec:Reg-cvrs}
In this section we define regular covers of Zariski structures and look 
into their structure, and in particular analyse the relations between fibres. 
The notion of regular cover we present here is compatible with the more 
general notion of cover of a first order structure given by Hrushovski 
in~\cite{Hr89} and repeated by Ahlbrandt and Ziegler in~\cite{AZ91}. The 
notion of regular cover presented here has some additional topological 
properties coming from Zariski structures, and as the main perspective 
is to establish a theory of specialisations; which requires the topological 
structure. One main difference is that, in regular cover structures the 
same group is acting on fibres. For general covers, there is more flexibility, 
different groups are allowed to act on different fibres.

\begin{definition}\label{def:z-cover-str}
Let $\mathcal C := (C,M,\pr)$ be a Zariski structure with two sorts $M$ and $C$ (called {\em base} and {\em cover} respectively) such that
\begin{enumerate}[label=(\roman*)]
\item There is a Zariski group $G$ in $\mathcal C$ acting morphically 
and freely on $C$ with Zariski continuous bijections. 
\item $M$ is interpretable in $C$ as a topological sort and $\pr:C\to M$ 
denotes the canonical quotient-map. It is an $\emptyset$-definable 
surjection. 
\item For each $m\in M$, the fibre $\pr^{-1}(m)$ is an orbit of an element 
in $G$, i.e. $\pr^{-1}(m) = G\cdot x$ for some $x\in C$.
\item The group $G$ is a Zariski group in $M$ (in particular 
$G$ is interpretable as a topological sort in $M$)
\end{enumerate}
Then $\mathcal C := (C,M,\pr)$ is said to be a \emph{regular cover (of $M$)}. 
\end{definition}

\begin{proposition}
Let $\mathcal C =(C, M,\pr)$ be a regular cover. Then
\begin{enumerate}[label=(\roman*)]
\item $G$ act transitively on each fibre of $\pr$ (hence it acts regularly);
\item the map $\pr: C\to M$ is a morphism;
\item $M$ is isomorphic to $C/G$; i.e. there is a bijective morphism 
$C/{G}\to M$.
\end{enumerate}
\end{proposition}
\begin{proof}
\begin{enumerate}[label=(\roman*)]
\item Immediate.
\item The map $\pr$ is the natural quotient map. Hence it is a morphism.
\item Since the action of the group $G$ is free, it is also proper by definition 
(see Definition~\ref{def:cont-prop-act}, (ii) and (iii)). Therefore the equivalence 
relation $E_G$ defined by 
\[c E_G d \mbox{ if and only if } \exists g\in G\, (g\cdot c = d)\] 
is closed. Hence $C/{E_G}$ is a topological sort. Moreover, $\pr$ is constant 
on the classes of $E_G$. By Definition~\ref{def:z-cover-str} (iii), $E_G$ is 
the pre-image of $=$ under $\pr$. Then $\pr/{E_G}:C/{G}\to M$ is 
an isomorphism.
\end{enumerate}
\end{proof}


\begin{definition}
Let $\bar{b} = (b_1,\ldots,b_n) \in C^n$, and $\pr(\bar{b}) = (m_1,\ldots,m_n) = 
\bar{m}$. Let $A\subset C$. We say $\bar{b}$ is \emph{strongly independent 
in fibres over $A$} if 
\[loc(\bar{b}/M\cup A)=\{ \bar{c}\in C^n: \pr\, \bar{c}=\bar{m}\}.\]
\end{definition}

\begin{lemma}\label{lem:ind-in-fibres}
Let $\mathcal C= (C,M,\pr)$ be a regular cover, and $m\subset M$. 
Suppose there is $b^\prime\in \pr^{-1}(m)$ strongly independent in 
fibres over $A$ for some $A\subseteq\mathcal C$. Then every $b\in 
\pr^{-1}(m)$ is independent in fibres over $A$ and is generic over $A
\cup M$.
\end{lemma}
\begin{proof}
If $b$ is dependent, it satisfies for some $a\in A$, and $m^\prime\in 
M$ a positive quantifier free formula $Q(m^\prime,a,b)$, with 
$\dim Q(m^\prime,a,y)<\dim \pr^{-1}(m)$. Then $b^\prime\in g\cdot 
Q(m^\prime,a, C)$ for some $g\in G(M)$, which gives a similar formula 
for $b^\prime$.
\end{proof}

\begin{lemma}\label{lem:ind-fibres-loc}
Suppose $\mathcal C\preceq \mathcal C^\prime$, where $\bar{b}^\prime\subset 
C^\prime\setminus C$, and  $|\bar{b}^\prime| = n$; also let $\bar{b}^\prime$ be 
strongly independent in fibres over $A\subseteq C$. Let $m^\prime:=\pr(b^\prime)$. 
Then 
\begin{enumerate}[label=(\roman*)]
\item the locus of $b^\prime$  over $A\cup M^\prime$ is of the form $S(y,m^\prime)$ 
where 
\[S(y,x)\equiv x=\pr(y)\] 
\item the locus of $b^\prime m^\prime$ over $A\cup M$ is of the form
$S(y,x)\ \& \ R(x,a)$ for some $a\subset M$, and $R(y,z)$ a positive quantifier free 
formula over $\emptyset$.
\end{enumerate}
\end{lemma}
\begin{proof}
\begin{enumerate}[label=(\roman*)]
\item Immediate by definition.
\item Let $Q(y,x,z)$ be a positive quantifier free formula such that $Q(y,x,c)$ is 
the locus of $b^\prime m^\prime$ over $A\cup M$, where $c\subset A\cup M$. 

Let $R(x,a)$ be the locus of $m^\prime$ over $A\cup M$, where $a\subset A\cup 
M$. Since $M$ is totally transcendental and stably embedded in $\mathcal C$, we 
may choose $a\subset M$. Since $M$ is a submodel of $M^\prime$, the locus 
$R(x,a)$ is irreducible. We may assume $a\subset c$ and $Q(y,x,z)\equiv Q(y,x,z)\ \&\ R(x,z)$.

For $m\in M^n$, let $S_m(y)$ be the formula $S(y,m)$; which is equivalent to $\pr(y) 
= m$. Let $S_{m}$ be the fibre in the respective model. Set 
\[R^0(M,c)=\{ m\in M^n: S_m\cap Q(C,m,c)\neq \emptyset\}\]
as the projection of $Q(C,m,c)\subset C^n\times M^n$ on $M^n$.

Note that $R^0(M,c)$ is a dense subset of $R(M,c),$ since $m'$ is a generic point 
in the $M^\prime$ versions of both.

We now consider the action of $g\in G^n$ on the set $C^n\times M^n$
\[(b,m) \mapsto (g\cdot b,m)\]

By our assumptions it is continuous, and thus $g\cdot Q(C,M,c)$ is closed.

By (i), $Q(C^\prime,m^\prime,c) = S_{m^\prime}$, and this is a generic fibre. Hence 
$G^n$ acts transitively on the generic fibre $Q(C^\prime,m^\prime,c)$.
So for any $g\in G(M)$,
\[g\cdot Q(C^\prime,m^\prime,c)=Q(C^\prime,m^\prime,c)\] 
Hence, in $\mathcal C$, for any $g\in G$ and any $m\in R^0(M,c)$, if $m$ is generic 
over $g,c$, then 
\[g\cdot Q(C,m,c)=Q(C,m,c)\]
 
By the addition formula, for a given $g\in G(M)$, 
\begin{align*}
&\dim g\cdot Q(C,M,c)\cap Q(C,M,c)\ge \\
 &\ge \dim G^n + \dim \{ m\in R^0(M,c): g\cdot Q(C,m,c)=Q(C,m,c)\}\ge\\ 
& \ge \dim Q(C,M,c)
\end{align*}
and hence, since $Q(C,M,c)$ is irreducible,  
\[g\cdot Q(C,M,c)=  Q(C,M,c)\]
for all $g\in G(M)$. This proves that every fibre $Q(C,m,c)$ is stable under the action 
of $G$, hence 
\[Q(C,m,c) = S_m(C) \mbox{ for all } m\in R^0(M,c)\]

Then $\pr^{-1} R^0(M,c)= Q(C,M,c)$ and $Q(C,M,c)$ is a closed $G$-invariant set. 
By definition, the image $R^0(M,c)$ is closed in the sort $M^n$. Clearly then $R^0(M,c) 
= R(M,a)$ and thus $S(y,x)\ \& \ R(x,a)\equiv Q(y,x,c)$.
\end{enumerate}
\end{proof}

\begin{corollary}\label{cor:loc-fibre}
If $Q(y,x,z)$ is the locus of $(b,m,c)$ over $\emptyset$ and $Q(C,m,c) = S_m$ 
then $Q(C,m^\prime,c^\prime)\neq \emptyset$ implies $Q(C,m^\prime,c^\prime) 
=S_{m^\prime}$ for any $m^\prime,c^\prime$.
\end{corollary}
\begin{proof}
By assumption $G^n$ acts transitively on the generic fibres $Q(C,m^\prime,c^\prime)$. 
So for any $g\in G(M)$ we have 
\[g\cdot Q(C,m^\prime,c^\prime) = Q(C,m^\prime,c^\prime)\]
Hence, in $\mathcal C$, for any $g\in G$ and for any generic $m^\prime,c^\prime$ 
satisfying $\exists y\, Q(y,x,z)$ we have 
\[g\cdot Q(C,m^\prime,c^\prime) = Q(C,m^\prime,c^\prime)\]

As all element of $\mathcal C$ are named by convention 
it is the prime model. In particular that elements of $M$ are named. Therefore, 
$m^\prime, c^\prime$ is generic over any $g\in G(M)$. Then, $g\cdot Q(C,M,C) 
\cap Q(C,M,C)$ contains the original generic element $bmc$.
Hence, since $Q(C,M,C)$ is irreducible,
\[g\cdot Q(C,M,C) = Q(C,M,C)\]
for all $g\in G(M)$. Which shows that every fibre $Q(C,m^\prime,c^\prime)$ 
is invariant under the action of $G(M)$. 
Hence 
\[Q(C,m^\prime,c^\prime) = S_{m^\prime} \mbox{ when } \exists y\, Q(y,m^\prime,c^\prime)\]
\end{proof}

\begin{corollary}\label{cor:loc-fibre-2}
Under assumptions of Corollary~\ref{cor:loc-fibre}, $\exists y\ Q(y,x,z)$ defines 
a closed set.
\end{corollary}
\begin{proof}
The topology on the sort $M\times M\times C$ can be defined from the sort 
$C\times M\times C$ by the equivalence relation $(y,x,z) \sim (g\cdot y,x,z)$ 
and the corresponding action of $G$,
\[C/{G}\times M\times C\simeq M\times M\times C\]
Corollary~\ref{cor:loc-fibre} together with Lemma~\ref{lem:ind-fibres-loc} proves 
that $Q$ defines a $G$-invariant closed subset of $C\times M\times C$. Hence
$Q/{G}$ is closed. This is homeomorphic to the set defined by $\exists y\ Q(y,x,z).$
\end{proof}

\begin{lemma}\label{lem:pr-complete-type}
Under assumptions of Lemma~\ref{lem:ind-fibres-loc}, let $b^\prime$ be strongly 
independent over $A\subseteq C$ and $m^\prime = \pr( b^\prime)$. Then
$S_{m^\prime}(y)$ defines a complete type (an atom of Boolean algebra) over 
$M^\prime\cup A$.
\end{lemma}
\begin{proof}
Suppose the formula $S(y,m^\prime)$ does not define a complete type over 
$A\cup M^\prime$.

Then there is a positive quantifier-free $Q(y,x,z)$ over $\emptyset$ 
such that for some $a\subset M^\prime\cup A$,
\[\vDash \exists y\ Q(y,m',a) \ \& \ S(y,m')\ \ \& \ \neg Q(b',m',a)\]
Let $b^{\prime\prime}$ satisfy $Q(b^{\prime\prime},m^\prime,a)$.
We may assume that $Q$ and $a$ is such that $Q(y,x,a)$ is the locus 
of $b^{\prime\prime} m^\prime$ over $M^\prime\cup A$. So $Q(C^\prime,m^\prime,a) 
\subset S_{m^\prime}$. 

Since $G(M^\prime)$ acts transitively on the fibres, $b^\prime=g\cdot 
b^{\prime\prime}$ for some $g\in G(M^\prime)$. Hence $g\cdot Q(C^\prime,m^\prime,a) 
\subset S_{m^\prime}$ and $g \cdot Q(C^\prime,m^\prime,a)$ is Zariski 
closed set defined over $M^\prime\cup A$. Since $b^\prime$ is strongly 
independent in fibres over $A$, it follows 
that $g\cdot Q(C^\prime,m^\prime,a)= S_{m^\prime}$. Hence 
$Q(C^\prime, m^\prime,a) = S_{m^\prime}$. A contradiction.
\end{proof}

\begin{lemma}\label{lem:loc-over-b}
Let $b\in C^n$ be strongly independent in fibres over $\emptyset$. Let 
$m^\prime\subset {M^\prime}^k$ and $Q(y,x^\prime,x,z) $ is a positive 
quantifier free formula and $mw\subset M$, $m=\pr(b)$, such that 
$Q(y,x^\prime,m,w)$ is the locus of $bm^\prime$ over $mw$. Then there 
is a positive quantifier free formula $R(x^\prime,x,z)$ over $\emptyset$ 
such that
\[ Q(y,x^\prime,x,z) \equiv R(x^\prime,x,z)\ \& \ x=\pr(y)\]
\end{lemma}
\begin{proof}
Let $R(x^\prime,x,z)$ be the formula $\exists  y\ Q(y,x^\prime,x,z)$.
We claim
\[Q(y,x^\prime,x,z)\equiv R(x^\prime,x,z)\ \& \ x=\pr(y).\]
Indeed, the implication from left to right is obvious. To see the inverse we 
need to prove that for any $(a^\prime,b,a,c)$ 
\[\vDash \forall y\, (Q(b,a^\prime,a,c)\ \& \ a=\pr\, y)\to Q(y,a^\prime,a,c)\]
But this formula immediately follows from Lemma~\ref{lem:pr-complete-type}.
%
This proves the claim.

In turn, the claim implies that the closed subset $Q$ of $C^n\times M^l$, 
some $l$, is saturated with respect to the equivalence relation on $C^n$ 
given by the action of group $G^n$. Moreover, the subset $R$ defined by 
$R(x^\prime,x,z)$ on the topological sort $M^l$ (see~\ref{def:z-cover-str}) 
is by definition of the topology closed. Hence $R(x^\prime,x,z)$ is positive 
quantifier free.
\end{proof}

We are making the following assumptions, which are to be valid for the rest 
of the paper:
\begin{itemize}
\item[-] $G$ is definably almost simple (i.e. proper definable normal 
subgroups are finite); and any definable normal subgroup $H\ideal G^k$ 
(for any $k\in\mathbb N$) is definable without parameters.
\item[-] The group $G$ is definable in $M$, so that $G\subset M^l$ is 
a definable subset for some $l$.
\item[-] The fibres $\pr^{-1}(m),$ for all $m\in M$, are atoms over $M$. 
\end{itemize}

Generally the definition of definably almost simple requires the group to be 
non-abelian, and (definably) connected. Our version does not require these 
conditions. In other words, we are allowing abelian groups, and groups which 
are not connected to be definably almost simple. In fact, if $G$ is definably 
almost simple and not definably connected then it must be finite. Indeed, 
in this case the connected component $G^0$ is finite since it is a definable 
normal subgroup. Moreover, it has finite index by definition. Hence $G$ is 
finite.

\begin{lemma}\label{lem:reducn-to-ADS-interm}
Let $b= b_1 b_2$, 
where $b_1=(b_{11},\ldots,b_{1n})\in C^{n}$ is a tuple and $b_2$ is a 
singleton such that $b_{11},\ldots,b_{1n}, b_{2}$, is not strongly independent 
in fibres over $\emptyset$. Assume also that the locus $\phi(b_1, C)$ of $b_2$ 
over $Mb_1$ is a proper subset of the fibre containing $b_2$.

Then there is a finite, $\emptyset$-definable $H\ideal G$  such 
that 
\[H\cdot b_2 = \phi(b_1, C)\]
and $\phi(b_1, C)$ is an atom over $ M\cup \{b_1\}$. Moreover,
for any $b^\prime_1\equiv_M b_1$ there is $b^\prime_2$ such that 
\[H\cdot  b^\prime_2 = \phi(b^\prime_1, C)\]
\end{lemma}
\begin{proof}
Note first that under the assumptions $\pr^{-1}(m_2)$ is an atom over $M$.

Let $\phi({b_1},C)$ be the locus of $b_2$ over $Mb_1$. Observe that $\phi({b_1},C)$ 
is an atom over $Mb_1$. If not we may assume that there is a proper $Mb_1$-closed 
subset $\psi({b_1},C)\subset \phi({b_1},C)$. There is a $g\in G$ such that 
$b_2\in g\cdot \psi({b_1},C)$. But $g\cdot \psi({b_1},C)$ is $Mb_1$-closed.  

Define a binary relation $E_{b_1}$ on $\pr^{-1}(m_2)$ as follows:
\[E_{b_1}(x,y) :\Leftrightarrow \exists g\in G\ (x\in g\cdot\phi({b_1},C) \ \& \ y\in g\cdot\phi({b_1},C))\]
$E_{b_1}$ is an $Mb_1$-definable equivalence relation on $\pr^{-1}(m_2)$. All equivalence 
classes are shifts of $\phi({b_1},C)$ by elements of $G$.

Reflexivity and symmetry of $E_{b_1}$ is obvious. We need the following claim 
to prove transitivity:

\begin{subclaim} Let $g_1, g_2\in G$, and assume $g_1\cdot \phi({b_1},C)\,\cap\, g_2\cdot 
\phi({b_1},C)\neq\emptyset$. Then $g_1\cdot \phi({b_1},C) = g_2\cdot \phi({b_1},C)$. 
\end{subclaim}
\begin{proof}
Suppose $y\in g_1\cdot \phi({b_1},C)\, \cap\, g_2\cdot \phi({b_1},C)$. Then 
there is an $h\in G$ such that $h\cdot y = b_2$. Then, $b_2\in hg_1\cdot 
\phi({b_1},C)\, \cap\, hg_2\cdot \phi_{b_1}(C)$. Since $\phi({b_1},C)$ is the 
locus of $b_2$ over $Mb_1$
\[hg_1\cdot \phi({b_1},C) = hg_2\cdot \phi({b_1},C) = \phi({b_1},C) = 
loc(b_2/Mb_1)\]
Then it follows that $g_1\cdot \phi({b_1},C) = g_2\cdot \phi({b_1},C)$. 
This proves the claim.
\end{proof}

\medskip

Next we will show that $E_{b_1}(y,z)$ is transitive. Let $E_{b_1}(x,y)$ and 
$E_{b_1}(y,z)$. Then there are $g_1,g_2\in G$ such that 
\begin{align*}
x\in g_1\cdot \phi_{b_1}(C) &\ \&\ y\in g_1\cdot \phi_{b_1}(C)\\
y\in g_2\cdot \phi_{b_1}(C) &\ \&\ z\in g_2\cdot \phi_{b_1}(C)\\
\end{align*} 
Observe that $g_1\cdot \phi_{b_1}(C)\,\cap\, g_2\cdot \phi_{b_1}(C)\neq 
\emptyset$, namely $y$ is in the intersection. Then by the previous claim 
\[g_1\cdot \phi_{b_1}(C) = g_2\cdot \phi_{b_1}(C)\]
Therefore $x\in g_2\cdot \phi_{b_1}(C)$. Hence $E_{b_1}(x,z)$.
\medskip

\begin{subclaim} 
$E_{b_1}$ is $M$-definable. 
\end{subclaim}
\begin{proof}
Let $b^\prime_1\equiv_M b_1$. Then $\phi(b^\prime_1,C)$ is the locus of 
some $b^\prime_2$ over $Mb^\prime_1$. Then by repeating the argument 
above one sees that $E_{b^\prime_1}$, defined in the same way, is an equivalence 
relation whose classes are shifts of $\phi(b^\prime_1,C)$; which is an atom 
over $Mb^\prime_1$. In particular, $b_2\in g\cdot \phi(b^\prime_1,C)$ for 
some $g\in G$. Since $b^\prime_1 = f\cdot b_1$ for some $f\in G^{n}$, we 
see that $g\cdot \phi(b^\prime_1,C)$ is $Mb_1$-definable and hence, $g\cdot 
\phi(b^\prime_1,C) = \phi(b_1,C)$. It follows, $E_{b_1}$ and $E_{b^\prime_1}$ 
have the same classes, and thus are equal. The claim follows.
\end{proof}
 
\medskip







 
Let $E:=E_{b_1}.$  
Define the subset $H\subset G$ as follows: 
\[g\in H :\Leftrightarrow g\cdot b_2 \in \phi({b_1},C)\] 
Then $H$ is actually a subgroup: First note that $g\cdot \phi({b_1},C) = \phi({b_1},C)$. 
Then it is immediate that product of two elements of $H$ is again in $H$.
It also follows from this observation that $H$ is closed under inversion; 
let $h\in H$ and $b_2^\prime := h^{-1}\cdot b_2$. Then $b^\prime_2\in h^{-1} 
\cdot\phi(b_1,C)$. By definition of $H$ we also know that $b_2\in h^{-1}\cdot 
\phi(b_1,C)$. Since $h^{-1}\cdot\phi(b_1,C)\,\cap\, \phi(b_1,C)\neq \emptyset$, 
we know that $h^{-1}\cdot\phi(b_1,C) = \phi(b_1,C)$. Hence $b_2^\prime = 
h^{-1}\cdot b_2 \in \phi(b_1,C)$. Moreover, $H\cdot b_2 = \phi({b_1},C)$.

\begin{subclaim} 
For any $c\in \pr^{-1}(m_2)$ we have $E(c,C)= H\cdot c$.
\end{subclaim}
\begin{proof}
The set
\[\{ c \in \pr^{-1}(m_2):\ E(c,C)=H\cdot c\}\]
is 
$M$-definable (as $M$ is stably embedded in $C$) and contains $b_2$. 
Hence, it must  be equal to $\pr^{-1}(m_2)$ since the latter is an atom. 
This proves the claim.
\end{proof} 

It follows that $E(x,y)\equiv \exists h\in H\ y=h\cdot x$, that is $E$ is 
induced by the action of the definable subgroup $H$. By Fact~\ref{fct:clsr-subgr-constr-mono}, 
$H$ is a closed subgroup. Then by the proof of~\ref{lem:z-topl-grp-quot}, 
$E$ is a Zariski closed subset of $\pr^{-1}({m}_2)\times \pr^{-1}({m}_2)$.

Also, $H$ is normal. Indeed, any $E$-class has the form $H\cdot c$, 
and for any $g\in G$ we have a class $H\cdot (g\cdot c)$. On the other 
hand, by definition of $E$, we know that the action of an element of 
$G$ takes an equivalence class to another equivalence class. Therefore 
$g\cdot (H\cdot c)  = g\cdot (c/{E})$ is an equivalence class. Moreover, 
we have $g\cdot (H\cdot c)= H\cdot (g\cdot c) $ which implies 
\[gHg^{-1}=H\]

Since $G$ is assumed to be definably almost simple, $H$ is finite. Also, 
$H$ is $\emptyset$-definable by assumption on $G$. 
\end{proof}
 
\begin{theorem}\label{thm:reducn-to-ADS}
Let $b= b_1 b_2$, where $b_1=(b_{11},\ldots,b_{1n})\in C^{n}$ is strongly 
independent in fibres over $U$ and $b_2=(b_{21},\ldots,b_{2k})$ is such 
that $b_{11},\ldots,b_{1n}, b_{2i}$, for each $i=1,\ldots,k$ is not strongly 
independent. Let $m_1=\pr(b_1)$, and $m_2=\pr(b_2)$.

Then there are $\emptyset$-definable normal finite subgroups $H_i\ideal G$ 
for $i=1,\ldots,k$, a $\emptyset$-definable subgroup $H\leq H_1\times\cdots
\times H_k$, 
and a positive quantifier free formula $\phi(y_1, y_2)$ over $UM$ such 
that 
\[\vDash \phi(b_1, b_2)\mbox{ and } H\cdot b_2 =  \phi(b_1, C^k)\]
and  $\phi(b_1, C^k)$ is an atom over $ M\cup \{b_1\}\cup U$. 
\end{theorem}
\begin{proof}
First, let us consider the case $U = \emptyset$. For any $b_{2i}$, its 
locus $\phi_i(b_1,C)$ over $Mb_1$ is a proper subset of $\pr^{-1}(m_{2i})$. 
For $i=1,\ldots,k$, let $H_i$ be the finite, $\emptyset$-definable normal 
subgroup of $G$ given by Lemma~\ref{lem:reducn-to-ADS-interm}. 
Let  $\phi(b_1,C^k)$ be the locus of $b_2$ over $Mb_1$. Then 
\[\phi(b_1,C^k)\subset (H_1\times\cdots\times H_k)\cdot b_2\]
 
Observe that, similarly to the proof of Lemma~\ref{lem:reducn-to-ADS-interm}, 
$\phi(b_1,C^k)$ is an atom over $Mb_1$. From a similar argument it 
follows that there is an $Mb_1$-definable subgroup $H<G^k$ such 
that $\phi(b_1,C^k) = H\cdot b_2$. In fact, $H<H_1\times\cdots\times 
H_k$. Therefore $H$ is finite. In particular, it is a finite subset of the 
prime model; hence $\emptyset$-definable as elements of the prime 
model are named. This proves the theorem in case $U = \emptyset$.

%
%

In the general case we may assume that $U = u\subset C$ is a finite 
subset, and thus $b^\prime_1 := ub_1$ would satisfy the assumptions we 
used for $b_1b_2$ in Lemma~\ref{lem:reducn-to-ADS-interm} and 
its application.
\end{proof}

\begin{corollary}\label{cor:cnct-fn-1}
Under assumptions and notation of Theorem~\ref{thm:reducn-to-ADS} 
there is a $U$-definable set $D_u$ and function $f_u: D_u\to C^k/{H}$, 
such that, given $b_1b_2$ there is a tuple $s_b\subset M$ 

 $$ s_bm_1m_2b_1\in D_u$$ 

and 
\[f_u( {s_b}, m_1,m_2,{b}_1) =\hat{b}_2\in \pr^{-1}( {m}_2)/{H}, \ \ \mbox{ where } \hat{b}_2=H\cdot  {b}_2
\]
\end{corollary}
\begin{proof}
$f_u(x, x_1,x_2,{y}_1)= z$ is given by the formula 
\[\exists  {y}_2 (\varphi_u(x, x_1,x_2,{y}_1,  {y}_2)\ \& \  {z}=H\cdot  {y}_2)\] 
where $\varphi_u(s_b,m_1,m_2,y_1,y_2)=\phi(y_1,y_2),$ with $\varphi$ as 
in~\ref{thm:reducn-to-ADS}, and $\varphi_u$ over $U$;
    and
$D_u$ is defined by 
\[(x, x_1,x_2,y_1) \in D_u \Leftrightarrow \exists  y_2 \ \{ \varphi_u(x, x_1,x_2,y_1,  y_2)\ \& \ \forall  y_2^\prime \ 
(\varphi_u(x,x_1,x_2, y_1,  y_2^\prime)\to  y_2^\prime\in H\cdot  y_2)\}\]  
\end{proof}

\subsection{Examples}

 1. Consider  a Zariski structure $\mathcal{C}=(C,M,\pr)$ which is a principal bundle $\pr: C\to M$ over a complex manifold $\mathbf{M}$ with the structure group $G.$ Since the fibration is locally trivial there is locally a  section of $\pr,$
$$c: M\to C$$
a regular map. Hence, for close enough $x_1,x_2\in M$ we can define 
$z_1=c(x_1)$ and $z_2=c(x_2).$ The infinitesimal version of the relation 
$\Phi$ between $x_1,x_2,z_1$ and $z_2$ is called a connection on the bundle.

A slightly more complex situation arises when we consider more general choice of elements $y_1$ and $y_2$ in fibres over $x_1$ and $x_2$ respectively. Now we take into account that $y_1=g_1\cdot z_1$ and $y_2=g_2\cdot z_2$ for some $g_1,g_2\in G.$ Hence the connection between $
x_1,x_2,y_1$ and $y_2$ requires parameters $g_1,g_2$ which in  \ref{cor:cnct-fn-1} correspond to $s_b$ (or the free variable $x$). Note that by construction 
 $g_1,g_2\in \dcl(x_1,x_2,y_1,y_2),$ that is $g_1,g_2$ are functions of $x_1,x_2,y_1,y_2.$ 

\medskip

2. This class of examples $\mathcal{C}=(C,M,\pr)$ belongs to the category of finite \'etale covers $\pr: C\to M$ of a smooth quasi-projective algebraic variety defined over a field $\mathrm{k}\subset \mathrm{F},$ where $\mathrm{F}$ is algebraically closed. The universe $M=\mathbf{M}(\mathrm{F})$ is the set of $\mathrm{F}$-points of the algebraic variety, and the structure on $\mathbf{M}$ is the standard Zariski structure associated with a quasi-projective algebraic variety. The formalism of \'etale covers assumes that the structure on the sort $C$ is definable in the field $\mathrm{F}$ with parameters in $\mathrm{k}$ but is not bi-interpretable with this structure. $\mathcal{C}$ satisfies the following conditions
\begin{itemize}
\item[(i)] $\dcl(M)\cap C=\emptyset.$ 

Moreover, the orbit of every $c\in C$ under $G=\mathrm{Aut}(\mathcal{C}/M)$ is the whole fibre $\pr^{-1}(\pr(c)).$ 
\item[(ii)] for any $c\in C,$ 
$$\dcl(\{ c\}\cup M)\supseteq C.$$
\end{itemize} 

Condition (i) implies that $\mathcal{C}$ is a regular cover with the finite group $G.$  

Condition (ii) implies that for any $b_1,b_2\in C$ there is $s_b\subset M,$ such that $b_2\in \dcl(s_b,b_1),$
that is $b_2=f(s_b,b_1)$ for some rational map $f$ over $\mathrm{k}.$

For technical convenience assume that the substructure $\mathbf{M}$ along with $M$ contains also a sort $\mathrm{F}$ for the field. Then we can choose $s_b\subset \mathrm{F}.$ Now we can use elimination of imaginaries in  algebraically closed fields to claim that $s_b\in \dcl(b_1,b_2).$

\medskip

The above examples  from algebraic geometry are generalised 
in the assumption (CC) below.  It will play an essential role in the 
next section, which is the core of this paper. In fact, some form of 
CC is necessary in order for the specialisations of $\mathcal{C}$ 
to have the nice model-theoretic properties proved in our main theorem. 

\subsection{Continuous Connections (CC) Assumption}\label{ssect:CC}

We now introduce the most important technical assumption.  

\begin{description}
\item[(CC)]\label{redD} 
Under assumptions and notation \ref{thm:reducn-to-ADS}
and \ref{cor:cnct-fn-1}:
\begin{itemize}[leftmargin=*]
\item[-] The restriction of $f_u(x,x_1,x_2,y_1)$ 
on a closed $u$-definable subset $D_u^0\subset D_u$ 
 containing $(s_b,m_1,m_2,b_1)$  is a morphism.
\item[-] There is a closed $u$-definable set $D_u^\dagger$ 
 containing $(m_1,m_2,b_1, \hat{b}_2)$ and a morphism
$$f_u^\dagger:D_u^\dagger\to M^{|s_b|}$$
 such that,
 if $(x,x_1,x_2,y_1) \in D_u^0$ and $f_u(x,x_1,x_2,y_1)=\hat{y}_2$, 
then $(x_1,x_2,y_1, \hat{y}_2) \in D_u^\dagger$ and  $f^\dagger(x_1,x_2,
y_1, \hat{y}_2)=x$.
\end{itemize}
\end{description}

{\bf Remark.} It is clear from the construction that $x_1\in \dcl(y_1)$ and 
$x_2\in \dcl(\hat{y}_2).$ We still keep the variables $x_1$ in $f_u$ and 
$x_1,x_2$ in $f_u^\dagger$ to make the notation graphical.  

\begin{lemma}\label{lem:rglr-good-models}
Let $\mathcal C = (C,M,\pr)$ be a regular cover satisfying the 
Continuous Connections assumption, then every model of 
$\mathrm{Th}(\mathcal C)$ is a regular cover and satisfies (CC).
\end{lemma}
\begin{proof}
Let $\mathcal D\models \mathrm{Th}(\mathcal C)$. First let us check that 
$\mathcal D$ is a regular cover. It is clear that $\mathcal D$ is 
a two sorted Zariski structure. We will write $\mathcal D= (D,N)$. The 
interpretation of the topological Zariski group $G$ gives a Zariski topological 
group $B$ interpretable in $N$. Also $B$ will act morphically and freely 
with Zariski automorphisms as all of these are first order properties. This 
is enough to see that $\mathcal D$ is a regular cover. 
By construction, $\mathcal D$ satisfies (CC).  
\end{proof}

\begin{remark}
The paper \cite{HZ96} presents the (historically first) example of a non-classical Zariski geometry, see section 10 of \cite{HZ96}.
Renaming $X^*$ of the example as $C$ and $X$ as $M$ we identify a cover structure
 $\mathcal C = (C,M,\pr).$ It is not difficult to see that $\mathcal{C}$ is a regular cover satisfying the 
Continuous Connections assumption (CC).

Many more Zariski cover structures satisfying (CC) can be found in  \cite{SSZ} and \cite{Zilber08}.
\end{remark}

\section{Specialisations of Regular Covers of Zariski Structures}\label{sec:spcl-cvrs}
In this section we will work with both languages; the Zariski language 
$\mathcal{L}$ of a (multi-sorted) Zariski structure, and $\mathcal{L}^{\pi}$, 
the language $\mathcal{L}$ expanded by a function $\pi$ which will 
be interpreted as a specialisation. 

\begin{lemma}\label{lem:univ-finiteness}
Let $\pi:\mathcal{C}\to \mathcal{C}_0$ be a specialisation, where 
$\mathcal C$ is a (multi-sorted) Zariski structure and $\mathcal{C}_0
\preceq \mathcal{C}$ .
Suppose that every $\aleph_0$-saturated model of the $\mathcal L^{\pi}$-theory 
$\mathrm{Th}(\mathcal C, \mathcal C_0, \pi)$ is $\aleph_0$-universal. 
Then every $\kappa$-saturated model of the theory is $\kappa$-universal.
\end{lemma}
\begin{proof}
Let $(\mathcal D,\mathcal D_0,\pi)$ be a $\kappa$-saturated model 
of the theory $\mathrm{Th}(\mathcal C, \mathcal C_0, \pi)$. So in 
particular it is $\aleph_0$-saturated. Hence by assumption $\pi$ is 
$\aleph_0$-universal. Let $\mathcal D^\prime \succeq \mathcal D$, 
$A\subseteq \mathcal D^\prime$ with $|A|<\kappa$ and $\pi^\prime: 
A\cup\mathcal D\to \mathcal D_0$ a specialisation extending $\pi$. 
Let $A_0:= A\cap\mathcal D$, and $A^\prime := A\setminus A_0$. Also, 
for any element of $A$, without loss of generality we may assume that 
its image under $\pi^\prime$ is in $A_0$. 

Let $a^\prime\subset A^\prime$ be an arbitrary finite tuple, and $B\subset 
A_0$ be an arbitrary finite subset. Let 
$\pi^\prime(a^\prime) =: a_0$. 
We need to show that the $\mathcal L^{\pi}$-type 
\[p_{a^\prime}(x/B) := \mathrm{tp}(a^\prime/B)\cup \{\pi(x) = a_0\}\]
is satisfiable in $\mathcal D$, where $\mathrm{tp}(a^\prime/B)$ is the 
$\mathcal L$-type. 

Since $a^\prime$ is finite and $\pi^\prime$ is $\aleph_0$-universal, 
there is an embedding $\sigma:a^\prime B\to \mathcal D$ over $B$ 
such that $\pi^\prime(a^\prime) = \pi(\sigma a^\prime)$. Set $a=\sigma(a^\prime)$. 
Then, $a\models \mathrm{tp}(a^\prime/B)$, and  $\pi(a) = a_0 = 
\pi^\prime(a^\prime)$.

By compactness, the type \[\displaystyle\bigcup_{\substack{a^\prime\subset 
A^\prime\\ B\subset A_0}} p_{a^\prime}(x/B)\] where $a^\prime$ ranges over all 
finite tuples in $A^\prime$ and $B$ ranges over all finite subsets of $A_0$ 
is realisable in some elementary extension of $(\mathcal D,\mathcal D_0,\pi)$. 
But since $(\mathcal D,\mathcal D_0,\pi)$ is $\kappa$-saturated, in fact there 
is a realisation of this type in $(\mathcal D,\mathcal D_0,\pi)$.
\end{proof}
 
In this section we will study specialisations of regular covers of Zariski
structures. Let $\mathcal C_0 = (C_0 ,M_0,\pr)$ be a regular cover (of 
the Zariski structure $M_0$), we will also often call $\pr:C_0\to M_0$ 
a regular cover. Let $\mathcal L$ denote the Zariski language for the 
two sorted Zariski structure $(C_0, M_0)$. Let $\mathcal C = (C, M, \pr) 
\succeq (C_0 ,M_0,\pr) = \mathcal C_0$, and let $\pi:\mathcal C \to 
\mathcal C_0$ be a maximal specialisation such that its restriction 
$\pi_M: M \to M_0$ is an $\aleph_0$-universal specialisation. 
For the rest of the paper  the assumptions made here on $(\mathcal C,
\mathcal C_0, \pi)$ are valid; most importantly, that $\pi$ is maximal 
and $\pi_M$ is $\aleph_0$-universal. 

By $\mathcal L^{\pi}$ we denote the language $\mathcal L$ expanded 
with a symbol $\pi$ which will be interpreted as the specialisation. 
In fact, we will consider $M_0$ along with the definable Zariski group 
$G$ as assumed previously. We will write $G(M_0)$ to indicate 
the realisation of $G$ in the Zariski structure $M_0$.  Therefore the 
elementary extension $M$ will be considered with the corresponding 
definable Zariski topological group $G(M)$. Often we will also consider 
$\pi_M$ as a specialisation $G(M)\to G(M_0)$ of the Zariski topological 
groups in the natural way.

In particular, if we take $\mathcal{C}$ to 
be a regular cover (of some Zariski structure $M$), then it suggests 
that we can restrict ourselves to study $\aleph_0$-universal specialisations 
of regular covers of Zariski structures. 

We can make a further reduction. It is enough to consider $\mathcal C^\prime
\succeq \mathcal C$, a finite subset $A\subset M\cup \Dom{\pi}$, and 
a finite tuple $b^\prime$ in $\mathcal C^\prime$, a specialisation $\pi^\prime: 
\mathcal C\cup \{b^\prime\} \to \mathcal C_0$ extending $\pi:\mathcal C
\to \mathcal C_0$. We aim to show that there is a $b$ in $\mathcal C$ such 
that 
\begin{equation} \label{eqn:univ-special-cover}
\pi(b)=\pi^\prime(b^\prime)\mbox{ and }\mathrm{tp}(b^\prime /A) = 
\mathrm{tp}(b/A)
\end{equation} 

By Theorem~\ref{thm:charac-univ-spc}, this will imply that $\pi$ is 
$\aleph_0$-universal. 
In fact, we show that if the regular cover, in addition to the assumptions 
made above, satisfies (CC) one can always 
find such a $b$ that satisfies~(\ref{eqn:univ-special-cover}). At the end 
of the section we present a characterisation of $\aleph_0$-universal 
specialisations of regular covers satisfying (CC). Which, in turn, allows 
us to weaken the assumption of Lemma~\ref{lem:univ-finiteness} for 
regular covers satisfying (CC).

\begin{proposition}\label{prop:cover-max-then-univ}
Let $\mathcal C_0\prec \mathcal C$ be a Zariski structure and its 
extension, let $\pi: \mathcal C\to \mathcal C_0$ be an $\aleph_0$-universal 
specialisation. Assume that $M$ is a sort in $\mathcal C$ and $\pr: 
C\to M$ is a regular cover satisfying (CC).

Let $\pi_M: M\to M_0$ be the restriction of $\pi$ to the substructure. 
Then $\pi_M$ is $\aleph_0$-universal.
\end{proposition}
\begin{proof}
Let $M^\prime \succ M$, and $n^\prime\subset M^\prime$ be a finite 
tuple. Let $\pi^\prime_M: M^\prime\to M_0$ be a specialisation extending 
$\pi_M: M\to M_0$ with $n^\prime\subset \Dom{\pi^\prime_M}$, and 
$\pi^\prime_M(n^\prime) = n_0$.
 
We claim that $\pi^\prime_M\cup \pi$ is a specialisation extending $\pi$. 
This will imply that the type of $n^\prime$ over a finite subset of $M$ is 
realised in $\mathcal C$, and so in $M$, by some $n$ such that $\pi(n) 
= n_0$. In other words, that $\pi_M$ is $\aleph_0$-universal.

Consider a positive quantifier free formula $Q(x^\prime,y)$ in $\mathcal C$ 
such that $Q(x^\prime,b)$, for some $b$ in $\Dom\pi$, is the locus of 
$n^\prime$ over $\Dom{\pi}$. Let $n := \pr(b)$. Then it follows $n\in 
\Dom\pi$.

Let $b = b_1b_2$ be the splitting of $b$ into $b_1$, which is maximal 
independent in fibres over $\emptyset$, and $b_2$ such that 
\[\hat{b}_2= f(s_{b},m_1,m_2,b_1), \ \ 
s_{b}= f^\dagger(m_1,m_2,b_1, \hat{b}_2)\]
where $m_1=\pr(b_1)$ and $m_2=\pr(b_2)$. Note that $m_1,m_2\in\Dom{\pi}$, 
since $b_1,b_2\in\Dom{\pi}$.

We have  $\hat{b}_2\in \Dom\pi$, since $b_2\in \Dom\pi$. Then $s_{b}\in 
\Dom\pi$, since $f^\dagger$ is a morphism. Now we replace $Q(x^\prime,y)$ 
with the formula
\[Q^*(x^\prime,x,w,y_1)\equiv \exists y_2\ Q(x^\prime,y_1y_2) \ \&\ 
[\hat{y}_2=f(w,x,y_1)]\]
where $[\hat{y}_2=f(w,x,y_1)]$ is the positive quantifier free formula (see~\ref{redD})
\[(w,x,y_1)\in D^0\ \&\ (x,y_1\hat y_2)\in D^\dagger\ \&\ \hat y_2 = f(w,x,y_1).\]
In particular $\hat y_2 = H\cdot y_2\in C^{|b_2|}/H$, where $H$ is the finite group 
given in relation to $b_2$ by Theorem~\ref{thm:reducn-to-ADS}.


We claim that $Q^*$ defines a Zariski closed set. Let $k=|x^\prime xw|$ and 
$t=|y|=|y_1y_2|$. Define the function
\[\mathrm{id}\times \mathrm{p}:M^k\times C^t\to M^k\times C^{|y_1|}\times C^{|y_2|}/{H} \]
by $\mathrm{id}\times \mathrm{p}(x^\prime xwy_1y_2) = (x^\prime xwy_1\hat y_2)$; 
i.e. it is identity on $M^k\times C^{|y_1|}$ and $\mathrm{p}:C^{|y_2|}\to C^{|y_2|}/H$ 
is the canonical quotient map. The image $M^k\times C^{|y_1|}\times 
C^{|y_2|}/{H}$ can be identified with the sort $(M^k\times C^t)/{H}$ 
by taking the action of $H$ on $M^k$ to be trivial. Since $H$ is finite, 
$(M^k\times C^t)/{H}$ is an orbifold; and $\mathrm{id}\times\mathrm{p}$ 
becomes the canonical quotient map. 

It follows that $Q^*(x^\prime,x,w,y_1)$ is the image (in the orbifold $(M^k\times C^t)/{H}$) 
under $\mathrm{id}\times\mathrm{p}$ of the closed set defined by 
\[Q(x^\prime,y_1y_2) \ \& \ [\hat{y}_2= f(w,x,y_1)]\]
Then by Lemma~\ref{lem:quot-fin-z-top-grp}, $Q^*(x^\prime,x,w,y_1)$ 
is closed. 

Next we claim that $Q^*(x^\prime,n,s_b,b_1)$ defines the locus 
of $n^\prime$ over $\Dom{\pi}$. It is enough to prove that 
\[Q^*(M^\prime,n,s_b,b_1)\subseteq Q(M^\prime,b)\]
Note that by the construction of formula $Q^*$,  
\[\vDash Q^*(n^{\prime\prime},n,s_b,b_1)\Rightarrow \ \ \vDash 
Q(n^{\prime\prime},b_1b^\prime_2)\mbox{ for some }b^\prime_2\in \hat{b}_2\]
where $\hat{b}_2 = H\cdot b_2$.
By Theorem~\ref{thm:reducn-to-ADS},  $H\cdot b_2$  is an atom over 
$M^\prime b_1$. Hence $\vDash Q(n^{\prime\prime},b_1b^\prime_2)$ if and 
only if $\vDash Q(n^{\prime\prime},b_1b_2)$. This completes the proof 
of the claim. 

Clearly, 
\[\vDash Q^*(n^\prime,n,s_b,b_1)\]
and $Q^*(x^\prime,x, w,y_1)$ satisfies assumptions of Lemma~\ref{lem:loc-over-b}. 
Hence, 
splitting $x=x_1x_2,$ $|x_1|=|y_1|$ (in correspondence with $n = m_1m_2$) we 
get by Lemma~\ref{lem:loc-over-b}
\[Q^*(x^\prime ,x_1x_2, w,y_1)\equiv R(x^\prime ,x_1,x_2, w)\ \& \ x_1=\pr(y_1)\] 
for some positive quantifier free formula $R$.

This implies that the locus of $n^\prime$ over $\Dom \pi$ is determined by 
the formula $R(x^\prime,m_1,m_2, s_b)$, where  $m_1,m_2,s_b\subset\Dom{\pi_M}$. 
Since $\pi^\prime_M(n^\prime) = n_0$
we get
\[\vDash R(n_0, \pi(m_1),\pi(m_2),\pi(s_b))\]
and hence 
\[Q(n_0, \pi(b))\] 
This proves that $\pi^\prime_M\cup \pi$ preserves positive quantifier free 
formulae over $\Dom\pi$ and hence is an extension of specialisation $\pi$.
\end{proof}

\begin{lemma}\label{lem:inf-cover-univ-base}
Let $\mathcal C^\prime\succeq\mathcal C$, $A\subset M\cup\Dom{\pi}$ be 
a finite subset, $b^\prime\subseteq \mathcal C^\prime$ be a finite tuple, 
and $\pi^\prime:\mathcal Cb^\prime\to\mathcal C_0$ be a 
specialisation extending $\pi$. Suppose $b^\prime \subset M^\prime$. Then 
we can find $b\subset M$ such that $\pi(b)=\pi^\prime(b^\prime)$ and 
$\mathrm{tp}(b^\prime /A) = \mathrm{tp}(b/A)$ (i.e. property~(\ref{eqn:univ-special-cover})
 holds).
\end{lemma}
\begin{proof}
Since $M^\prime$ is totally transcendental and stably embedded in $\mathcal C^\prime$, 
there is a finite $A_{M^\prime} \subset M^\prime$, such that $\mathrm{tp}(b^\prime 
/A_{M^\prime} ) \vdash \mathrm{tp}(b^\prime /A)$. Since the restriction $\pi^\prime_M: 
\{b^\prime\}\cup M \to M_0$ of $\pi^\prime$ is a specialisation extending $\pi_M$, 
and since $\pi_M$ is $\aleph_0$-universal, there is $b\subset M$ such that $\pi(b) 
= \pi^\prime(b^\prime)\mbox{ and }\mathrm{tp}(b^\prime/A_{M^\prime}) = 
\mathrm{tp}(b/A_{M^\prime})$. Property~\ref{eqn:univ-special-cover} follows.
\end{proof}

From here on we will assume that $M$ is $\aleph_0$ saturated.

\begin{lemma}\label{lem:ext-spcl-same-fibre}
Let $\mathcal C^\prime\succeq\mathcal C$ and $b^\prime\subseteq \mathcal C^\prime$ 
be a finite tuple. Let $\pi^\prime:\mathcal Cb^\prime\to\mathcal C_0$ be a 
specialisation extending $\pi$, and $A\subset M\cup\Dom{\pi}$ be a finite 
subset. Suppose that $\pr(b^\prime) = \pr(a)$ for some $a\subset C\cap
\Dom{\pi}$. Then there is $b\subset \mathcal C$ such that $\pi(b) = 
\pi^\prime(b^\prime)$ and $\mathrm{tp}(b^\prime /A) = 
\mathrm{tp}(b/A)$ (i.e. property~\ref{eqn:univ-special-cover} holds).
\end{lemma}
\begin{proof}
Since $b^\prime$ and $a$ are in the same fibre, $b^\prime = g^\prime 
\cdot a$ for some unique $g^\prime \in G(M^\prime)$. Since the group 
action is free, we may assume $g^\prime$ is in the domain of $\pi^\prime_M$. 
If not, using Lemma~\ref{lem:ext-spcl-via-morph}, one can extend 
$\pi^\prime_M$ to a specialisation $M^\prime\to M_0$ which is defined 
on $g^\prime$. With the abuse of notation we will again denote this extension 
by $\pi^\prime_M: M^\prime\to M_0$, and also write $g^\prime\in\Dom{
\pi^\prime_M}$. Let $a_0 := \pi(a)$, $b_0 = \pi^\prime (b^\prime)$, and 
$g_0 = \pi_M^\prime(g^\prime)$. 

It follows that $\pi_M^\prime$ is an extension of the specialisation 
$\pi_M$. 
Since $\pi_M$ is $\aleph_0$-universal there is 
$g \in G(M)$ such that $\mathrm{tp}(g/Aa) = \mathrm{tp}(g^\prime/Aa)$ 
and $\pi_M(g) = g_0$.

Let $b := g \cdot a$. By this definition, $\mathrm{tp}(b^\prime/Aa)$ 
and $\mathrm{tp}(b/Aa)$ are determined by $\mathrm{tp}(g^\prime/Aa)$ 
and $\mathrm{tp}(g/Aa)$ respectively. So we have $\mathrm{tp}(b^\prime/Aa) 
= \mathrm{tp}(b/Aa)$. Finally, again by freeness of the group action, 
$\pi(b)$ is defined. Then $\pi(b) = \pi(g) \cdot \pi(a) = g_0 \cdot a_0$, 
that is $\pi(b) = b_0.$ Hence our choice of $b$ satisfies property~\ref{eqn:univ-special-cover}. 
\end{proof}

\begin{remark}
Below, in the proofs of Lemma~\ref{lem:ext-spcl-to-gnrc}, Lemma~\ref{lem:spc-ext-base-to-fbr} 
and Proposition~\ref{prop:cov-spcl-main} we make a case distinction 
between $G$ being infinite and finite. This may look strange at first. As 
explained in Remark~\ref{rmk:non-conn-finite-grps}, the main difference 
here is between $G$ being connected and not connected. 
The only non-connected case in our setting is when $G$ is finite. By 
definition, when $G$ is infinite it is irreducible; as $G$ is a group, that 
means it is connected.

In addition, we keep the case distinction between finite and infinite to emphasize 
that these cases also correspond to the situation where the fibres of the cover are 
infinite or finite. 
\end{remark}

\begin{lemma}\label{lem:ext-spcl-to-gnrc}
Suppose $G$ is infinite. Let 
$m\in M^n\cap \Dom{\pi}$, with $m_0 = \pi(m)$, 
and $b_0\in \pr^{-1}(m_0)\cap C^n_0$. 
Let $\mathcal C^\prime\succeq \mathcal C$. Then, for any $b^\prime 
\in {C^\prime}^n\cap \pr^{-1}(m)$ generic over $\Dom{\pi}$ there is a 
specialisation $\pi^\prime:\mathcal Cb^\prime\to \mathcal C_0$ extending 
$\pi$ such that $\pi^\prime(b^\prime) = b_0$.
\end{lemma}
\begin{proof}
Note that since $G$ is connected, so is $G^n$; and hence $\pr^{-1}(m)$ 
is irreducible. Let $b^\prime\in {C^\prime}^n$ be an element of $\pr^{-1}(m)$ 
which is  generic over $C$.   


Consider a positive quantifier free formula $Q(y,x,z)$ over $\emptyset$ 
and a tuple $c$ from $\Dom{\pi}$ such that $\vDash Q(b^\prime,m,c)$.
We may assume without loss of generality that $Q$ defines the locus of 
$b^\prime,m,c$ over $\emptyset$. By the genericness assumption, 
$Q(y,m,c)\equiv \pr(y) = m$.

Hence the assumptions of Corollary~\ref{cor:loc-fibre} are satisfied. Then 
by Corollary~\ref{cor:loc-fibre-2} $\exists y\,  Q(y,x,z)$ defines a closed set.  

It follows 
\[Q(\mathcal C_0,\pi(m),\pi(c))\neq \emptyset \mbox{ and } Q(\mathcal C_0,\pi(m),\pi(c))=S_{m_0}\]
Hence $\vDash Q(b_0,m_0,\pi(c))$ for any such $Q(y,x,c)$.  
 

Set $\pi^\prime(b^\prime)=b_0$. By construction $\pi^\prime$ is a specialisation 
extending $\pi$.
\end{proof}

\begin{lemma}\label{lem:spc-ext-base-to-fbr}
Let $\mathcal C_0\preceq\mathcal C$ be a regular cover
and its extension. 
Suppose $m\in M^n\cap \Dom\pi$. Then there is $b\in 
\pr^{-1}(m)\cap \Dom\pi$.
\end{lemma}
\begin{proof}
It is enough to prove the statement for each coordinate of $m$, so we can 
assume $n=1$. We write $\pr^{-1}(m)$ in $C$ as $C_m$.

When $G$ is finite the statement follows from Lemma~\ref{lem:ext-spcl-to-sort-fin-grp}.
Assume $G$ is infinite and let $\mathcal C^\prime\succeq \mathcal C$ 
such that $\mathcal C^\prime$ is $|\Dom{\pi}|^+$-saturated. 
We will consider the following two cases.
\begin{description}[leftmargin=0cm]
\item[Case 1.] There is $b^\prime\in C^\prime_m$ independent in fibres over 
$\Dom{\pi}$. Then, by Lemma~\ref{lem:ind-in-fibres}, every $b\in C_m$ is 
independent in fibres and is generic over $\Dom{\pi}$. By Lemma~\ref{lem:ext-spcl-to-gnrc}
there is an extension $\pi^\prime$ such that $b\in \Dom{\pi^\prime}$. But 
$\pi$ is maximal and so $\pi^\prime(b) = \pi(b)$.
\item[Case 2.] Any $b^\prime \in C^\prime_m$ is dependent in fibres over 
$\Dom \pi$. Choose $b^\prime\in C^\prime_m$ generic over $\Dom{\pi}$. 
By Lemma~\ref{lem:ext-spcl-to-gnrc} 
there is an extension $\pi^\prime$, and $b^\prime\in\Dom {\pi^\prime}$.

Since $b^\prime$ is a singleton, according to notations of the 
paragraph~\ref{redD}, we have   
$b^\prime_1 = \emptyset$ and $b^\prime = b^\prime_2$. Further, there is 
$s_{b^\prime}$ in $M^\prime$ and $u\in \Dom\pi$ 
such that
\[f_u(s_{b^\prime},m)=\hat{b}^\prime,\ \ f^\dagger_u(m,\hat{b}^\prime) 
= s_{b^\prime}\]
The second equality together with the fact that $f^\dagger$ is a morphism 
implies that $\pi^\prime$ can be extended to $s_{b^\prime}$. With an 
abuse of notation we will write $s_{b^\prime}\in\Dom{\pi^\prime}$. Consider 
the restriction of $\pi^\prime$ to $M^\prime$ and denote it by $\pi_M^\prime: 
M^\prime\to M_0$. 

By Lemma~\ref{lem:inf-cover-univ-base}, there is a $s_b\in M$ such 
that $\mathrm{tp}(s_{b^\prime}/mu) = \mathrm{tp}(s_b/mu)$ 
and $\pi^\prime(s_{b^\prime}) = \pi(s_b)$. 
Set $\hat{b}= f_u(s_b,m)$. Then $\hat{b}$ is an element of 
the topological sort $C/{H}$ where $H < G$ is finite (see Theorem~\ref{thm:reducn-to-ADS}), 
and by continuity of $f_u$, we have $\hat b\in\Dom{\pi}$. Hence 
by Lemma~\ref{lem:ext-spcl-to-sort-fin-grp} applied to the sort $C/{H}$, 
there is a $b\in\Dom{\pi}$ such that $\hat{b} = H\cdot b$. 
\end{description}
\end{proof}

\begin{corollary}\label{cor:spc-same-fib}
Suppose $m\in M^n\cap \Dom \pi$, and $A\subset C$ 
with $|A|<\aleph_0$. Suppose there is $b^\prime\in \pr^{-1}(m)$ 
in $C^\prime$ and $\pi^\prime$ extends $\pi$ so that $b^\prime\in 
\Dom{\pi^\prime}$. Then there is $b\in \pr^{-1}(m)\cap \Dom \pi$ such 
that $\pi(b)=\pi^\prime(b^\prime)$ and $\mathrm{tp}(b^\prime/A) = 
\mathrm{tp}(b/A)$. In other words property~\ref{eqn:univ-special-cover} 
holds.
\end{corollary}
\begin{proof}
Follows from Lemma~\ref{lem:spc-ext-base-to-fbr} and Lemma~\ref{lem:ext-spcl-same-fibre}.
\end{proof}

Next we will consider property~\ref{eqn:univ-special-cover} in more 
detail. 
Let us recall that we are considering $\mathcal C^\prime\succeq \mathcal C$, 
a finite subset 
\begin{equation}
A\subset M\cup\Dom{\pi},
\end{equation}\label{eqn:parameter-set-spcl-cover}
 a finite tuple $b^\prime$, and 
a specialisation $\pi^\prime:\mathcal C\cup\{b^\prime\}\to \mathcal C_0$ 
extending $\pi:\mathcal C\to \mathcal C_0$.
Suppose $b^\prime$ satisfies property~\ref{eqn:univ-special-cover}. 
I.e. there exist a tuple $b\in\mathcal C$ such that 
\[\pi(b) = \pi^\prime(b^\prime) \mbox{ and } \mathrm{tp}(b^\prime/A) 
= \mathrm{tp}(b/A)\]
Split $b^\prime = b_1^\prime b_2^\prime \in {C^\prime}^{n+k}$ so that 
$b_1^\prime\in {C^\prime}^n$ is maximal strongly independent in fibres 
over $\Dom\pi$ and $b_2^\prime$ is the rest. 

Then 
by Continuous Connections, 
\begin{equation}\label{eqn:fmb}
\hat{b}_2^\prime= f_{u}(s_{b^\prime},m^\prime_1,m^\prime_2,
b^\prime_1)\mbox{ and }s_{b^\prime} = f^\dagger_u(m^\prime_1, 
m^\prime_2,b^\prime_1,\hat{b}^\prime_2)
\end{equation}
for the morphism $f_u$ over a finite $u\subset \Dom\pi$, 
for some tuple $s_{b^\prime}$ in $M^\prime$, where $m^\prime_1 
= \pr(b^\prime_1)$, and $m^\prime_2=\pr(b^\prime_2)$. And $\hat{b}^\prime_2 
= H\cdot b^\prime_2$ is an element of ${C^\prime}^k/{H}$ with $H$ an 
$\emptyset$-definable finite subgroup of $G^n$ (see Theorem~\ref{thm:reducn-to-ADS} 
and Corollary~\ref{cor:cnct-fn-1}).

\begin{proposition}\label{prop:cov-spcl-main}
Let $\mathcal C^\prime\succeq \mathcal C$.
Assume $b^\prime\in {C^\prime}^n$ is strongly independent in fibres 
over $\Dom\pi$, $m^\prime = \pr(b^\prime)$, $s_{b^\prime} \subset M^\prime$ 
and $\pi^\prime:\mathcal C^\prime \to \mathcal C_0$ a specialisation, 
defined on $s_{b^\prime} m^\prime b^\prime$, extending $\pi$. Let 
$A\subset M\cup \Dom\pi$ be of cardinality less than $\aleph_0$. Then 
there is $s_b m b \subset \Dom\pi$ in $\mathcal C$ such that
\[\pi(s_b m b) = \pi^\prime(s_{b^\prime} m^\prime b^\prime)\mbox{ and }
\mathrm{tp}(s_b m b/A)=\mathrm{tp}(s_{b^\prime} m^\prime b^\prime/A)\]
\end{proposition}
\begin{proof}
Since $\pi_M$ is $\aleph_0$-universal, there is $s_b m\subset M\cap 
\Dom\pi$ such that $\mathrm{tp}(s_{b^\prime} m^\prime/A) = 
\mathrm{tp}(s_b m/A)$ and $\pi^\prime(s_{b^\prime}
m^\prime) = \pi(s_b m^\prime) =:  s_0 m_0$.  
\begin{description}[leftmargin=0cm]
\item[Case 1] If $G$ is finite. Then, by Lemma~\ref{lem:ext-spcl-to-sort-fin-grp}, 
we get $\pr^{-1}(m)\subset \Dom\pi$.

Since the sizes of fibres are equal and $\pi$ preserves the discrete 
Zariski topology on the fibres, 
\[\pi (\pr^{-1}(m)) = \pr^{-1}(m_0)\] 
and $\pi$ is a bijection on $\pr^{-1}(m)$.
In particular, there is $b\in \pr^{-1}(m)$ such that $\pi(b)=b_0$.  

We claim that $\mathrm{tp}(s_{b^\prime} m^\prime b^\prime/A) 
= \mathrm{tp}(s_b mb/A)$. This follows from the fact that $\mathrm{tp}
(s_{b^\prime} m^\prime/A)$ expresses that $S(y,m^\prime)$ is 
an atom over $s_{b^\prime} A$, and hence, by equality of types, 
$S(y,m)$ is an atom over $s_b A$. (Here $S(y,x)\equiv x=\pr(y)$, 
see Lemma~\ref{lem:ind-fibres-loc}).
\item[Case 2] If $G$ is infinite. Then By Lemma~\ref{lem:ext-spcl-to-gnrc}, 
for any $b^{\prime\prime}\in \pr^{-1}(m)$ that is generic over $C$ there is 
an extension $\pi^{\prime\prime}$ of $\pi$ such that $\pi^{\prime\prime}(b^{\prime\prime}) 
= b_0$. Note that $\mathrm{tp}(b^{\prime}/Am^\prime)$ is generic.

By Corollary~\ref{cor:spc-same-fib}, there exists $b\in \pr^{-1}(m)
\cap \Dom\pi$ satisfying the generic type over $Am$ in $\pr^{-1}(m)$ 
with $\pi(b) = b_0$. The type $\mathrm{tp}(b^\prime/Am^\prime)$ 
is generic in $\pr^{-1}(m^\prime)$ by assumption on $b^\prime$. 
Also $\mathrm{tp}(m^\prime/A) = \mathrm{tp}(m/A)$ by the above 
choice, hence $\mathrm{tp}(b^\prime/A) = \mathrm{tp}(b/A)$.
\end{description}
\end{proof}

\begin{corollary}\label{cor:main-rslt-1}
Let $A\subset M\cup\Dom{\pi}$ be of cardinality less than $\aleph_0$.
Assume that $\mathcal C$ satisfies the Continuous Connections assumption. 
Then property~(\ref{eqn:univ-special-cover}) holds for any $b^\prime \in 
{C^\prime}^n\cap(\Dom{\pi^\prime})^n$, where $\mathcal C^\prime\succeq 
\mathcal C$ and $\pi^\prime:\mathcal C^\prime\to \mathcal C_0$ is a specialisation 
extending $\pi$.
\end{corollary}
\begin{proof}
Write $b^\prime = b^\prime_1 b^\prime_2$ where $b^\prime_1$ is maximal 
strongly independent in fibres over $A$, and $b^\prime_2$ is the rest.
We need to find $b\in\mathcal C$ such that 
\[\pi(b) = \pi^\prime(b^\prime) \mbox{ and } \mathrm{tp}(b/A)=\mathrm{tp}(b^\prime/A)\]  

Let $\pr(b^\prime_1) =: m^\prime_1, \pr(b^\prime_2) =: m^\prime_2$. 
By Continuous Connections, 
\[\hat{b}_2^\prime= f_{u}(s_{b^\prime},m^\prime_1,m^\prime_2,
b^\prime_1)\]
for a morphism $f_u$ over a finite $u$, and $s_{b^\prime}\subset M^\prime$ 
is some finite tuple. We may assume that $m^\prime_1, m^\prime_2, s_{b^\prime}
\in \Dom{\pi^\prime}$.

Now, taking $Au$ instead of $A$ in Proposition~\ref{prop:cov-spcl-main} 
and repeating the same argument for $s_{b^\prime}m^\prime_1 m^\prime_2 
b^\prime_1$ gives
\begin{align*}
\mathrm{tp}(s_b, m_1,m_2,b_1/Au) = \mathrm{tp}(s_{b^\prime}, 
m^\prime_1, m^\prime_2, b_1^\prime/Au)
\end{align*}  
and 
\begin{align}\label{piAu}
\pi(s_b, m_1,m_2,b_1) = \pi^\prime(s_{b^\prime}, m^\prime_1,m^\prime_2,
b^\prime_1)
\end{align} 

Set $\hat{b}_2 := f_u(s_b, m_1,m_2,b_1)$. It follows from 
the equality between types, that 
\[\mathrm{tp}(b_1\hat{b}_2/Au) = \mathrm{tp}(b^\prime_1\hat{b}^\prime_2/Au)\]
By applying Lemma~\ref{lem:ext-spcl-via-morph} to $f_u^\dagger$, 
and (\ref{piAu}) we get $\pi(\hat{b}_2) = \pi^\prime(\hat{b}^\prime_2)$.

The set $\hat{b}_2=H\cdot b_2$ is an atom over $AMb_1$ by 
Theorem~\ref{thm:reducn-to-ADS} . So any choice of such a $b_2$ 
satisfies 
\[ \mathrm{tp}(b_1{b}_2/Au) = \mathrm{tp}(b^\prime_1,{b}^\prime_2/Au)\]
By Lemma~\ref{lem:ext-spcl-to-sort-fin-grp} we can choose 
$b_2$ such that $\pi(b_2) = \pi^\prime(b^\prime_2)$. 
\end{proof}

All of the analysis of specialisations of Zariski cover structures carried 
in this section yields the following theorem. Which gives a characterisation 
of $\aleph_0$-universal specialisations of Zariski cover structures satisfying 
the Continuous Connections assumption. 

\begin{proposition}\label{thm:cover-equiv-univ}
Let $\mathcal C_0 = (C_0,M_0,\pr),$  $\mathcal C = (C,M,\pr)$  be 
two regular cover structures satisfying Continuous Connections assumption,  $\mathcal C_0\preceq \mathcal C$
 and 
$\mathcal C$ is an $\aleph_0$-saturated extension. 
Let $\pi:\mathcal C\to \mathcal C_0$ be a specialisation.
Then the following are equivalent:
\begin{enumerate}[label=(\roman*)]
\item $\pi$ is  $\aleph_0$-universal; 
\item  the restriction $\pi_M: M \to M_0$
is $\aleph_0$-universal and $\pi: \mathcal C\to \mathcal C_0$ is a maximal specialisation; 
\item the restriction $\pi_M: M\to M_0$ 
is $\aleph_0$-universal, 
and the following sentences hold
\begin{align}\label{e2}
\forall m\in \Dom{\pi_M}\ \exists c\in \Dom\pi\  \pr(c)=m.
\end{align} 
\begin{align}\label{e2.1}
\forall c\in \Dom\pi\ \forall g\in\Dom{\pi_G} \ g\cdot c\in\Dom{\pi}.
\end{align}
\end{enumerate}
\end{proposition}
\begin{proof}
First we note that $\mathcal C$ and $M$ satisfy the assumptions 
of Lemma~\ref{lem:prime-&-mnml}. Hence we may use 
Theorem~\ref{thm:charac-univ-spc}(ii) as the criterion for universality. 
So, in order to prove ((ii)$\Rightarrow$(i)) we need to satisfy 
property~(\ref{eqn:univ-special-cover}) for $A\subset M\cup \Dom\pi$. 
This is Corollary~\ref{cor:main-rslt-1}. 

((i)$\Rightarrow$(ii)) follows from Proposition~\ref{prop:cover-max-then-univ}, 
since any $\aleph_0$-universal specialisation is maximal. 

((ii)$\Rightarrow$(iii)) Sentence~(\ref{e2})  follows from 
Lemma~\ref{lem:spc-ext-base-to-fbr}. Sentence~(\ref{e2.1}) 
follows from the maximality of $\pi$ and Lemma~\ref{lem:ext-spcl-via-morph}.

((iii)$\Rightarrow$(ii)) It is enough to show that $\pi$ is maximal. 
First, observe that $\pi_M$ is maximal, since it is $\aleph_0$-universal. 
Now suppose $\pi^\prime: \mathcal C\to \mathcal C_0$ is a specialisation 
extending $\pi$ with $c^\prime\in\Dom{\pi^\prime}\setminus\Dom{\pi}$. 
Let $m:=\pr(c^\prime)$, since $\pr$ is a morphism, one can extend 
$\pi_M$ to $m^\prime$ (by Lemma~\ref{lem:ext-spcl-via-morph}). 
But $\pi_M$ is maximal, hence $m\in\Dom{\pi_M}$.

By~(\ref{e2}), there is a $c\in\Dom{\pi}$ with $\pr(c) = m$. There 
is a unique $g\in G$ such that $g\cdot c^\prime = c$. By freeness 
of the action, we may extend $\pi^\prime_M$ to $g$. As before 
since $\pi_M$ is maximal we actually have $g\in \Dom{\pi_M}$. 
In particular, $g\in\Dom{\pi_G}$. Hence by~(\ref{e2.1}), $c^\prime 
\in\Dom{\pi}$.  
\end{proof}

\begin{proposition}
Let $\mathcal C_0 = (C_0,M_0,\pr)\preceq \mathcal C = (C,M,\pr)$ be 
two regular covers. Let $\pi:\mathcal C\to \mathcal C_0$ be a specialisation, 
such that it restriction $\pi_M:M\to M_0$ to $M$ is maximal. Suppose 
that sentences~(\ref{e2}) and~(\ref{e2.1}) of Theorem~\ref{thm:cover-equiv-univ} 
(iii) hold. Then for any $m\in\Dom{\pi_M}$, there is a $c\in\pr^{-1}(m)$ 
such that
\[\Dom{\pi}\cap \pr^{-1}(m) = (G\cap\Dom{\pi})\cdot c.\] 
\end{proposition}
\begin{proof}
By~(\ref{e2}) there is an element $c\in\Dom{\pi}\cap \pr^{-1}(m)$. 
Then for any $g\in G\cap \Dom{\pi}$, we have $g\cdot c\in \Dom{\pi}
\cap \pr^{-1}(m)$ by~(\ref{e2.1}). We remark here that $G\cap 
\Dom{\pi}$ is always non-empty (as $G(M_0)\subseteq \Dom{\pi_M}$). 

Now suppose $d\in\Dom{\pi}\cap \pr^{-1}(m)$ with $d\neq c$. 
Then there is a $g\in G$ such that $g\cdot c = d$. Hence we 
may extend $\pi$ to $g$ since the action is free. By construction, 
this will be an extension of $\pi_M$ (recall Lemma~\ref{lem:ext-spcl-via-morph}). 
But $\pi_M$ is maximal. Then $g\in\Dom{\pi_M}\subset\Dom{\pi}$. 
This proves $\Dom{\pi}\cap \pr^{-1}(m) = (G\cap \Dom{\pi})\cdot c$. 
\end{proof}

\begin{remark}
Under the assumption $\pi:\mathcal C\to \mathcal C_0$ is maximal 
and its restriction $\pi_M:M\to M_0$ is $\aleph_0$-universal,
sentence~(\ref{e2}) implies that $\pi(S_m(C)) = 
S_{\pi(m)}(C_0)$, since along with $c\in \Dom\pi$ we have $G(M_0)\cdot c
\subset \Dom\pi \cap S_m(C)$.
\end{remark}

\begin{proposition}
Under assumptions of Proposition~\ref{thm:cover-equiv-univ}, consider 
the structure $(\mathcal C,\mathcal C_0,\pi)$ in the language of 
specialisations and its substructure $\pi_M: M\to M_0$. Suppose 
$\pi:\mathcal C\to \mathcal C_0$ is maximal, and that 
any $\kappa$-saturated 
model of the theory of specialisations of $\pi_M: M\to M_0$ is 
$\kappa$-universal. 

Then any $\kappa$-saturated model of the theory of specialisations of 
$\pi:\mathcal C\to \mathcal C_0$ is $\kappa$-universal. 
\end{proposition}
\begin{proof}
Consider the theory of specialisation $\pi: \mathcal C\to \mathcal C_0$.
Pick a $\kappa$-saturated model of the theory. By Lemma~\ref{lem:rglr-good-models} 
we may assume it is  $\pi:\mathcal C\to \mathcal C_0$. Then, by our 
assumptions, $\pi_M: M\to M_0$ is $\kappa$-universal. Now by Theorem~\ref{thm:cover-equiv-univ} 
we get that $\pi$ is $\aleph_0$-universal. Lemma~\ref{lem:univ-finiteness} 
completes the proof.  
\end{proof}

\section{The Theory of Universal Specialisations for Regular Covers of Zariski Structures}\label{sec:theory}
We present a theory, $\mathrm{Th}(\mathcal C)^{\pi}$, of universal 
specialisations of regular covers of Zariski structures satisfying (CC) 
in the language $\mathcal L^{\pi}$ of specialisations.

Let $\mathcal C_0$ be a regular cover satisfying (CC), and $\mathcal C
\succeq\mathcal C_0$ be an elementary extension. Let $\pi:\mathcal C\to 
\mathcal C_0$ be a specialisation such that its restriction $\pi_M:M\to M_0$ 
is an $\aleph_0$-universal specialisation. In particular, this implies $\pi_M$ 
is non-trivial. Hence $\pi$ is non-trivial. 

\begin{proposition}
Assume that the theory of specialisation $\mathrm{Th}(M,M_0,\pi)$ admits quantifier elimination in $\mathcal{L}^{\pi}$. 
Also assume that $\pi:M\to M_0$ is $\aleph_0$-universal. Then every $\kappa$-saturated 
model of $\mathrm{Th}(M,M_0,\pi)$ is $\kappa$-universal.
\end{proposition}
\begin{proof}
An argument similar to the proof of Lemma~\ref{lem:strd-model-th} is enough 
to prove the claim for $\kappa = \aleph_0$. Then one can get the general case 
by Lemma~\ref{lem:univ-finiteness}.
%
\end{proof}

We now describe the theory  $\mathrm{Th}(\mathcal C)^{\pi}$ of specialisation of the cover structure. It consists of the following 
axioms:
\begin{itemize}
\item[T1.] The complete theory of the pair of Zariski cover structures $\mathcal C_0 
\subsetneq \mathcal C,$   $\mathcal C_0 =(C_0, M_0,\pr),$ $\mathcal C =(C, M,\pr)$.
\item[T2.] The specialisation axioms stating 
that $\pi:\mathcal C\to \mathcal C_0$, 
  for any positive quantifier free 
$\mathcal L$-formula $Q(x)$ that 
\[\forall c\in \Dom \pi  \ Q(c)\to Q(c^\pi).\]
\item[T3.] The restriction $\pi_M: M\to M_0$ of $\pi$ to the base sort $M$ satisfies the
complete theory of  universal specialisation.
\item[T4.] The sentence
\[\forall m\in \Dom \pi\ \forall a\in \pr^{-1}(\pi(m))\ \exists c\in \pr^{-1}(m)\ \pi(c)=a\]
\item[T5.] The sentence 
\[\forall c\in \Dom\pi\ \forall g\in\Dom{\pi_G} \ g\cdot c\in\Dom{\pi}.\]
\item[T6.] The sentence 
\[\forall d\ (d\in\Dom{\pi}\implies \pr(d)\in\Dom{\pi}) \]
\end{itemize}

We will denote the models of $\mathrm{Th}(\mathcal C)^\pi$ with 
gothic letters $\mathfrak{B}, \mathfrak{C}, \mathfrak{D}, ...$ etc. 
More precisely, for a model $\mathfrak D\models \mathrm{Th}(\mathcal C)^\pi$, 
when we wish to emphasise the underlying structure, we will write 
$\mathfrak D = (\mathcal D,\mathcal D_0,\pi)$ where the sorts 
$\mathcal D$ and $\mathcal D_0$ are regular covers of Zariski 
structures with $\mathcal D_0\preceq \mathcal D$.
 
Whenever we consider any model $\mathfrak D = (\mathcal D, 
\mathcal D_0,\pi)$ of $\mathrm{Th}(\mathcal C)^\pi$, we will 
consider it together with all its topological sorts over $\emptyset$. 
Which in particular means that we will consider $\pi$ as a specialisation 
extending to the sorts in the natural way, via the corresponding 
quotient maps.

Let $(\mathcal A,\mathcal A_0,\pi)\subset (\mathcal D, \mathcal D_0,\pi)$ 
be a substructure. For a topological sort $T$ in $D$ over $\emptyset$ we 
will consider its relativisation to $A$. Recall $T = W/{E}$ for some 
$\emptyset$-definable $W\subseteq D^n$, relativise $W$ to $A$ as 
\[W_A := W\cap A^n = W(A)\] 
where $W(A)$ is the realisation of $W$ in $A$. Then the relativisation 
of $T$ to $A$ is 
\[T_A := W_A/{E}.\]  

\begin{lemma}\label{lem:qe-srts-spcl}
Let $\mathfrak D = (\mathcal D, \mathcal D_0, \pi)$ and $\mathfrak B 
= (\mathcal B, \mathcal B_0, \pi)$ be two models of $\mathrm{Th}(\mathcal C)^{\pi}$. 
Let $(\mathcal A,\mathcal A_0,\pi)\subset \mathfrak D$ be a substructure 
and $i:(\mathcal A,\mathcal A_0,\pi)\to \mathfrak B$ be a partial 
embedding. Then, the embedding $i$ can be extended to the (relativised) 
topological sorts $T_A$.
\end{lemma}
\begin{proof}
Let $T_A = W_A/{E}$ be a relativised topological sort. As $E$ is a 
closed equivalence relation it is preserved under $i$. As an $\mathcal L$-embedding, 
$i$ is actually an elementary embedding. Therefore $i(a)\in W(B)$ 
for any $a\in W_A$. Hence define 
\begin{align*}
i/{E} : T_A &\to T(B) \\
a/{E} &\mapsto i(a)/{E}.
\end{align*}
It is clear that $i/{E}$ preserves closed sets of $T_A$. It is also clear 
that $i\cup i/{E}: \mathcal A\cup T_A\to \mathcal B\cup T(B)$ is an embedding 
of Zariski structures. 
For the sake of notation we write $i$ instead of $i\cup i/{E}$. Next, 
we show that $i: \mathcal A\cup T_A\to \mathcal B\cup T(B)$ preserves 
the specialisation. 

Let $a\in\Dom{\pi}\cap T_A$. Then there is an $\alpha\in\Dom{\pi}
\cap W_A$ with $\alpha/{E} = a$ such that $\pi(a) = \pi(\alpha)/{E}$. 
Then $i(\pi(a)) = i(\pi(\alpha)/{E}) = \pi(i(\alpha))/{E} = \pi(a)$.
\end{proof}

\begin{theorem}\label{thm:QE}
Assume that the specialisation theory $\mathrm{Th}(M,M_0,\pi_M)$ of the base sort eliminates quantifiers. Then $\mathrm{Th}(\mathcal C)^{\pi},$  the theory of specialisation of the cover structure satisfying (CC),
admits quantifier elimination and  
is complete. 
\end{theorem}
\begin{proof}
Let $\mathfrak D = (\mathcal D, \mathcal D_0, \pi)$ and $\mathfrak B 
= (\mathcal B, \mathcal B_0, \pi)$ be two models of $\mathrm{Th}(\mathcal C)^{\pi}$ 
such that $\mathfrak B$ is $|\mathcal D|^+$-saturated. Let 
$(\mathcal A, \mathcal A_0,\pi) \subseteq \mathfrak D$ be a 
substructure and $i: (\mathcal A, \mathcal A_0, \pi)\to \mathfrak B$ 
be a partial embedding.

We will extend $i$ to an embedding $\mathcal A\cup \mathcal D_0\to \mathfrak B$. 
Observe that, as an $\mathcal L^{\pi}$-embedding, $i$ maps 
$\mathcal A_0 :=\mathcal A \cap \mathcal D_0$ to $\mathcal B_0$. 
Consider $i_{|_{A_0}}: \mathcal A_0\to \mathcal B_0$ as an embedding 
in the Zariski language $\mathcal L$ (without $\pi$). By assumption 
$\mathcal B_0$ is $|\mathcal D_0|^+$-saturated with respect to the 
language $\mathcal L$.  By quantifier elimination in this language, 
$i_{|_{A_0}}$ extends to an $\mathcal L$-embedding $j_0: \mathcal D_0
\to \mathcal B_0$. Since the specialisations are identity on both 
structures $\mathcal D_0$ and $\mathcal B_0$, the embedding $j_0$ 
is actually an $\mathcal L^{\pi}$-embedding.

Then it is immediate that $i\cup j_0: \mathcal A\cup D_0\to \mathfrak B$ 
is an $\mathcal L^{\pi}$-embedding. By abuse of notation we will write 
$i:\mathcal A\cup D_0\to \mathfrak B$ instead of $i\cup j_0$ and simply 
assume that $\mathcal A$ contains $\mathcal D_0$.

%

Let $N:=\pr(D),$  $d\in D\setminus A$ 
and $n := \pr(d)$.

Since $\mathrm{Th}(M,M_0,\pi_M)$ 
admits quantifier elimination, we can extend $i$ to  $\{ n\}\cup\pr(A)\to \pr(B)$ 
as an elementary monomorphism. Let $m:=i(n).$

We need to find an element $b\in B$, with $\pr(b)=m$ such that $i: d\mapsto b$ 
is an extension of the given embedding.

If $d\in \Dom \pi$ then $n\in \Dom \pi$ and hence $m\in \Dom \pi$. It is enough 
to find $b$ such that $nd$ and $mb$ satisfy the same $\mathcal{L}$-type over 
$A$ and so that $b\in \Dom \pi$ with $\pi(b) = \pi(i(d)) = i(\pi(d))$.
 
In  the easy case when $G$ is finite any $b\in \pr^{-1}(m)$ is in $\Dom\pi$ by 
Lemma~\ref{lem:ext-spcl-to-sort-fin-grp}. Hence, in particular there is a $b$ 
which satisfies the latter condition. 

Continuing with the finite case, if $d\notin \Dom \pi$ then by Lemma~\ref{lem:ext-spcl-to-sort-fin-grp} 
we have $n\notin \Dom \pi$ and hence $m\notin \Dom \pi.$ Again any $b\in \pr^{-1}(m)$ 
satisfying the $\mathcal{L}$-type satisfies the condition. 

Thus we are done when $G$ is finite.  

Now we assume that $G$ is infinite and hence, by our assumptions, connected.
There remain two cases:

\begin{description}[leftmargin=0cm]
\item[Case 1] $\pr^{-1}(n)$ is an $\mathcal L$-atom over $NA$. 
In this case one can extend the restriction $i:\pr(A)\to \pr(B)$ to an
embedding $i:\{ n\}\cup\pr(A)\to \pr(B)$ by the same argument as above.
Let $i(n) = m$.

If $d\in\Dom{\pi}$, put $d_0 :=\pi(d)$, and $b_0 := i(d_0)$. We will have $n\in\Dom{\pi}$, and hence $m\in\Dom{\pi}$. Then 
$i(\pi(n)) = \pi(m) = \pr(b_0)$. By axiom T4, there exists $b\in 
\pr^{-1}(m)$ such that $\pi(b) = b_0$.  Hence, $dn\mathcal A\to 
\mathcal B$ with $dn\mapsto bm$ is an extension of $i$. Since 
$\pr^{-1}(n)$ is an $\mathcal L$-atom over $NA$  the extension is an $\mathcal L^\pi$-embedding.

If $d\not\in\Dom{\pi}$, then it is enough to find a $b\in\pr^{-1}(m)
\setminus\Dom{\pi}$. By axiom T5 
\begin{align}\label{T5}
G(N)\not\subseteq \Dom{\pi} \mbox{ implies } \pr^{-1}(n)\not\subseteq\Dom{\pi}.
\end{align}
As this is a direct consequence of the axioms, the analogous 
statement with $m$ and $M$ instead of $n$ and $N$ will also 
be true. 

The left hand side of~(\ref{T5}) can be expressed as the type 
\[p(x) = \{x\in G\}\cup \{\pi(x)\neq t : t\in D_0\}\]
over $D_0$. The corresponding type in $(B,B_0,\pi)$ over $i(D_0)$ 
is 
\[q(x) = \{x\in G\}\cup\{\pi(x)\neq a : a\in i(D_0)\}\]
Since $(B,B_0,\pi)$ is $|\mathcal D|^+$-saturated type $q(x)$ is 
realised in this model. Then, by~(\ref{T5}) this means $\pr^{-1}(m)
\not\subseteq\Dom{\pi}$. Now, pick an element $b\in\pr^{-1}(m)
\setminus\Dom{\pi}$, and extend $i$ to $dnA\to B$ by sending 
$dn\mapsto bm$. As before this is enough to see the extension 
is an $\mathcal L^\pi$-embedding.

\item[Case 2] $\pr^{-1}(n)$ is not and $\mathcal L$-atom over 
$NA$. So, in particular $\bar e d$ is not strongly independent in fibres 
over $\emptyset$ for some finite tuple $\bar e\subset A$. 

\begin{subclaim} 
We can assume $\bar e$ to be strongly independent over $\emptyset.$ In other words,
there is $\bar{e}_1\subset \bar{e}$ such that $\bar{e}_1$ is strongly independent in fibres over $\emptyset$
and $\bar e_1 d$ is not strongly independent in fibres 
over $\emptyset.$ 
\end{subclaim}
\begin{proof}
Note that by Lemma~\ref{lem:reducn-to-ADS-interm}, 
$\dim(d/\bar{e}\cup N)=0$  (equivalently, the Morley rank).

Let $\bar{e}=\bar{e}_1\bar{e}_2$ where $\bar{e}_1$ is maximal 
strongly independent in $\bar{e}$, and we assume $\bar{e}_2\neq \emptyset$. 
Then $\bar{e}=\bar{e}_1\bar{e}_2$ satisfies the assumptions of 
Theorem~\ref{thm:reducn-to-ADS}, and hence there is a formula 
$\varphi(\bar{y}_1,\bar{y}_2)$ over $N$ such that $\varphi(\bar{e}_1,\bar{y}_2)$ 
is an atom over $N\cup \bar{e}_1$ realised by finitely many tuples, 
in particular by $\bar{e}_2$.  It follows that $\dim(\bar{e}_2/N\cup \bar{e}_1)=0$ .
 Hence  $\dim(d/\bar{e}_1\cup N)=0<\dim \pr^{-1}(n)$. Hence 
 $d\bar{e}_1$  is not strongly independent in fibres over $\emptyset$. 
\end{proof}

Then by 
Theorem~\ref{thm:reducn-to-ADS}, and Continuous Connections 
assumption, there is a finite subgroup $H\leq G$, a $\emptyset$-definable 
closed subsets $D^0$ and $D^\dagger$, and $\emptyset$-definable 
morphisms $f:D^0\to C/{H}$ and $f^\dagger:D^\dagger\to M$ 
such that 
\[\hat d = f(s_d,\bar k, n, \bar e) \mbox{ and } s_d = f^\dagger(\bar k, n, \bar e,\hat d)\] 
where $\bar k =\pr(\bar e)$, 
$\hat d = H\cdot d$, and $s_d$ is some tuple in $N$.

Since $\mathrm{Th}(M,M_0,\pi_M)$ admits quantifier elimination, 
one can extend $i$ to an embedding $ns_d\pr(A)\to Z$. Which in 
turn extends $i$ to an $\mathcal L^\pi$-embedding $ns_dA\to B$. 
Say $i(s_d) = r, i(\bar k) = \bar l, i(n) = m$ and, $i(\bar e) = \bar a$. 
Since $D^0$ is a $\emptyset$-definable closed set, we have $\models 
D^0(r,\bar l, m, \bar a)$. Put $\hat b := f(r,\bar l,m, \bar a)$. 

Now we extend $i$ to $(A\hat d)/{H}$ in the topological sort 
$C/{H}$ via the canonical quotient map as in Lemma~\ref{lem:qe-srts-spcl}. 
We claim that (the extension of) $i$ maps $\hat d\in A\hat d/{H}$ to $\hat b$. 
In Lemma~\ref{lem:qe-srts-spcl} we established that $i:A\cup (A\hat d/{H})\to 
\mathfrak B$ is an $\mathcal L$-embedding, hence preserves the closed 
subsets of $A^n\times (A\hat d/{H})^m$. In particular it will preserve the graphs 
of $f$ and $f^\dagger$, which are closed as $f$ is a morphism. Hence 
we get 
\begin{align*}
i(\hat d) &= i(f(s_d,\bar k, n, \bar e)) = f(r,\bar l,m, \bar a) =\hat b\\
i(s_d) &= i(f^\dagger(\hat d, \bar k, n,\bar e)) = f(\hat b, \bar l, m, \bar a) = r.
\end{align*}
By construction $\hat b = H\cdot b$ for some $b$.
 


If $d\in\Dom{\pi}$, put $\pi(d) = d_0$, and $i(d_0) = b_0$. Then 
we also see that $\pi(\hat d) = \hat d_0$ and $i(\hat d_0) = \hat b_0$. 
We claim that $\hat b\in\Dom{\pi}$ and $\pi(\hat b) = \hat b_0$. If not, 
extend $\pi:B/H\to B_0/H$ to $\pi^\prime :B/H\to B_0/H$ by defining 
$\pi^\prime(\hat b) = \hat b_0$. 

\begin{subclaim}
$\pi^\prime :B/H\to B_0/H$ is a specialisation.
\end{subclaim}
\begin{proof}
Let $S\subset (B/H)^n$ be a closed subset. Assume $S(\hat b, c)$, 
and $\hat b, c\in\Dom{\pi^\prime}$. By construction $S = T/H$ for 
some closed set $T$ of $B$. So, $T(b, \gamma)$ for all $b\in \hat b$ 
and $\gamma\in H\cdot c$. Observe also that $T/G$ is a closed set 
of $M$. Then, $T/G(m,v)$, where $\pr(b) = m$ and $\pr(\gamma) = v$. 
Then $T/G(\pi_M(m),\pi_M(v))$. Hence, $T(b_0, \gamma_0)$. Then 
$S(\hat b_0, \hat c_0)$.  
\end{proof}

\begin{subclaim}
$\pi^\prime\cup \pi : B\cup B/H\to B_0\cup B_0/H$ is a specialisation.
\end{subclaim}
\begin{proof}
Immediate.
\end{proof}

\begin{subclaim}
$\pi^\prime\cup \pi : B\cup B/H\to B_0\cup B_0/H$ can be extended to 
a specialisation $\pi^{\prime\prime} : B\cup B/H\to B_0\cup B_0/H$ such 
that $H\cdot b\subset \Dom{\pi^{\prime\prime}}$.
\end{subclaim}
\begin{proof}
Observe that $B\to B/H$ (and so $B_0\to B_0/H$) is an orbifold with t
the structure group $H$. Now, an argument similar to the one used 
in the proof of Lemma~\ref{lem:ext-spcl-to-sort-fin-grp}, one can show 
that $\pi\cup\pi^\prime$ can be extended to $\pi^{\prime\prime}$ such 
that $H\cdot b\subset \Dom{\pi^{\prime\prime}}$ the whole fibre $H\cdot b$. 
\end{proof}

Which in particular means that $\pi:B\to B_0$ can be extended to 
$H\cdot b$. But $\pi$ is $\kappa$-universal, hence maximal. So 
$H\cdot b\subset \Dom{\pi}$. Then, in particular there is a $b^\prime 
\in H\cdot b$ with $\pi(b^\prime) = b_0$. Now we can extend $\sigma$ 
to $d$ by $\sigma(d) = b^\prime$.

If $d\not\in\Dom{\pi}$, then we claim that $\hat b\cap\Dom{\pi} = 
\emptyset$. First we will consider the case $\pr(d) = n\in\Dom{\pi}$. 
By axiom T4 there is a $d^\prime\in\pr^{-1}(n)\cap \Dom{\pi}$. Also 
there is a $g\in G$ such that $g\cdot d^\prime = d$. By axiom T5, 
$g\not\in \Dom{\pi}$. We may assume that $d^\prime$ and $g$ are 
in $\mathcal A$.

Consider the type
\begin{align*}
\mathrm{p}(x) = \{x\in G\}\cup \{\pi(x)\neq a : a\in D_0\}\cup\{x\cdot d^\prime = d\}.
\end{align*}
By the above paragraph it is realised in $\mathfrak{ D}$.

Consider the corresponding type 
\begin{align*}
\mathrm{q}(x) = \{x\in G\}\cup \{\pi(x)\neq a : a\in i(D_0)\}\cup\{x\cdot b^\prime = b\}
\end{align*}
in $\mathfrak{B}$.

By stable embeddedness we may assume that $\mathrm{q}(x)$ has 
only parameters in $M$. By quantifier elimination of  $\mathrm{Th}(M,M_0,\pi_M)$, 
it is also consistent. By saturation of $\mathfrak{B}$, we see that 
$\mathrm{q}(x)$ is realised in $M$. It now follows that $b\not\in\Dom{\pi}$. 
In fact it also follows that $\hat b\cap\Dom{\pi} = \emptyset$. Hence 
the extension of $i$ to $dn\mathcal A\to \mathfrak B$ given by 
$dn\mapsto bm$ is an $\mathcal L^\pi$-embedding.  

If $n\not\in\Dom{\pi}$, then $m\not\in\Dom{\pi}$. Then it follows 
that $\pr^{-1}(m)\cap\Dom{\pi} =\emptyset$. Then the argument 
follows as above. 

\end{description}
This establishes the quantifier elimination. Next we show that 
$\mathrm{Th}(\mathcal C)^\pi$ is complete. Consider $(\mathcal C_0,
\mathcal C_0,\id)$ where $\id:\mathcal C_0\to\mathcal C_0$ 
is the identity map. Pick an element $m\in M\setminus M_0$, 
and add it to $C_0$ together with the whole fibre $\pr^{-1}(m)$. 
Then $(\mathcal C_0\pr^{-1}(m)m, \mathcal C_0,\pi)$ where 
$\pi:\mathcal C_0\pr^{-1}(m)m\to \mathcal C_0$ is the restriction 
of $\pi:\mathcal C\to \mathcal C_0$. Hence $(\mathcal C_0\pr^{-1}(m)m, 
\mathcal C_0,\pi)$ is a prime substructure of the theory 
$\mathrm{Th}(\mathcal C)^\pi$. Together with quantifier elimination, 
this implies that $\mathrm{Th}(\mathcal C)^\pi$ is 
complete.
\end{proof}
\begin{corollary} Let $\mathcal{C}=(C,M,\pr)$ be a regular cover satisfying (CC) and assume $M$ carries the standard Zariski structure associated with an algebraic variety over an algebraically closed field. Then $\mathrm{Th}(\mathcal C)^{\pi}$ is complete and admits elimination of quantifiers.   
\end{corollary}

\begin{lemma}\label{lem:strd-model-th}
Any $\aleph_0$-saturated model of the theory $\mathrm{Th}(\mathcal C)^{\pi}$ 
is $\aleph_0$-universal.
\end{lemma}
\begin{proof}
Let $(\mathcal D,\mathcal D_0,\pi)$ be an $\aleph_0$-saturated model 
of $\mathrm{Th}(\mathcal C)^{\pi}$. Let $\mathcal D^\prime\succeq \mathcal D$ 
be an elementary extension with respect to the Zariski language $\mathcal L$ 
(without $\pi$). Let $A\subset \mathcal D^\prime$ be a finite set, and 
$\pi^\prime: A\cup\mathcal D\to \mathcal D_0$ be a specialisation 
extending $\pi$. Without loss of generality we assume that $\pi^\prime(a)
\in A$ for each $a\in A\cap\Dom{\pi}$.

Enumerate $A\setminus \mathcal D$ as $\{a_i: i < k \}$ where $k\leq |A|$. 
Let $A_0 = A\cap \mathcal D$ and $\sigma_0 = \id :A_0\to \mathcal D$. 
Define $A_j :=A_0\cup\{a_i: i < j\}$. 

Let $i = j+1$, and assume that a partial elementary embedding $\sigma_j:A_j
\to \mathcal D$ over $A_0$ with $\pi^\prime(a) = \pi(\sigma(a))$ for 
all $a\in A_j\cap\Dom{\pi^\prime}$ is constructed. Consider the $\mathcal L^\pi$ 
type 
\[\mathrm{p}(x) = \{\varphi(x,\sigma_j(\bar a)) : \mathcal D^\prime\models\varphi(b_i,\bar a), \mbox{ and } \bar a\in A_j\}.\]

Since $A_j$ is finite this type is realised in $(\mathcal D,\mathcal D_0,\pi)$ 
by some $b_i$, since the structure is $\aleph_0$-saturated. Define $\sigma_i(a_i) = b_i$. 
Observe that $\pi^\prime(a_i)\in A_0$, say $\pi^\prime(a_i) = \alpha$. 
Then $\pi(x) = \alpha \in\mathrm{p}(x)$. Hence $\pi(\sigma_i(a_i)) = 
\pi(b_i) = \alpha$. It follows that $\sigma_i$ is a a partial elementary 
map with the desired property. Hence we are done by induction, proving 
$\pi:\mathcal D\to \mathcal D_0$ is $\aleph_0$-universal. 
\end{proof}

\begin{corollary}
Any $\kappa$-saturated model of $\mathrm{Th}(\mathcal C)^{\pi}$ 
is $\kappa$-universal.
\end{corollary}
\begin{proof}
Follows from Lemma~\ref{lem:strd-model-th} and Lemma~\ref{lem:univ-finiteness}. 
\end{proof}
\section{The 1996-Example of a non-algebraic Zariski geometry as a Regular Cover with CC}\label{sect:HZ-ex}

We consider here the example of a non-classical Zariski structure introduced in~\cite{HZ96}, section 10 and show that it
satisfies (CC).

Let $M$ be a one dimensional Zariski geometry, and $B\leq ZAut(M)$ be 
a group acting freely on $M$ by Zariski automorphisms. Let $B^*$ be a 
group extension of $B$ with a finite kernel $G$:
\[1\to G\to B^*\to B\to 1\]

Let $C$ be a set such that $B^*$ acts freely on it in the same way $B$ 
acts on $M$. I.e. the number of $B^*$ orbits is the same as number of 
$G$ orbits. 
 
For $c,d\in C$ define $c\equiv d$ if and only if there is an $h\in G$ 
such that $h\cdot c = d$. Clearly $\equiv$ is an equivalence relation invariant 
under the action of $B^*$. Moreover $\equiv$ is closed (hence definable). 
So $G$ acts trivially on the quotient $C/{\equiv}$, and hence the action 
of $B^*$ will give an action of $B$ on this quotient. By construction 
$C/{\equiv}$ is isomorphic to $M$ as $B$-sets, say via a map $\pr$. 
This map $\pr$ can be extended naturally to $\pr:C_0\to M_0$ so that 
$\pr(g^* c) = p(g^*)\pr(c)$ for all $g^*\in B^*$, where $p:B^*\to B$ is 
the group homomorphism. 

Then $C_0$ is made into a Zariski geometry by defining pull-backs of 
the closed sets of $M_0^n$ via $\pr$ and graphs of elements of $G^*$ 
as basic relations (of an associated language) and declaring the Boolean 
combinations of these predicates as closed sets. 

Let $\bar b = (b_1,\ldots, b_n)\in C^n$. Then $\bar b$ is strongly independent 
in fibres over $\emptyset$ if and only if $b_i$ is not in the orbit of $b_j$ for 
any $i,j$.  

Now let $\bar b\in C^n$, write $\bar b = b_1b_2$ where $b_1$ is maximal 
strongly independent in fibres over $\emptyset$ and $b_2$ is the rest. Let 
$\phi(b_1,C^k)$ be the locus of $b_2$ over $Mb_1$. By Corollary~\ref{sIndC}, 
there is a $\emptyset$-definable function $f:D\to C^k/H$ where $H<G^k$, 
$\{s_bm_1m_2\}\times \pr^{-1}(m_1)\subseteq D$, ${m}_1= \pr({b}_1), 
m_2 = \pr(b_2)$ and $s_b\subset M$ is a parameter. The function $f$ is 
defined as 
\[\exists  {y}_2 (\phi(x, x_1,x_2,{y}_1,  {y}_2)\ \& \  {z}=H\cdot  {y}_2)\]
and $D$ is defined as 
\[(x, x_1,x_2,y_1) \in D \Leftrightarrow \exists  y_2 \ \{ \phi(x, x_1,x_2,y_1,  y_2)\ \& \ \forall  y_2^\prime \ 
(\phi(x,x_1,x_2, y_1,  y_2^\prime)\to  y_2^\prime\in H\cdot  y_2)\}\] 

In the next paragraph we will see that the parameter $s_b$ is actually 
unnecessary. Later in the analysis, we will even conclude that $m_1,m_2$ 
are also not necessary as parameters. 

We claim that the locus $\phi(b_1,C^k)$ is the singleton $\{b_2\}$. 
Write $b_1 = (b_{11},\ldots,b_{1n})$ and $b_2 = (b_{21},\ldots, b_{2k})$. 
Since $b_1b_{2i}$, for all $i=1,\ldots,k$, are dependent in fibres, there is 
a $b_{1j}\in b_1$ and a $g_{ij}\in B^*$ such that $g_{ij}\cdot b_{1j} = b_{2i}$. 
The pair $g_{ij}$ and $b_{1j}$ uniquely determine $b_{2i}$. Therefore, 
$b_2$ is in the definable closure of $b_1$. Moreover, as the graphs of 
elements of $B^*$ are closed sets, we see that $\{b_2\}$ is a $b_1$-closed 
set (i.e. it is a closed set only using $b_1$ as a parameter). Hence the locus 
$\phi(b_1,C^k)$ of $b_2$ is the singleton $\{b_2\}$. Hence the only 
subgroup $H$ of $(B^*)^k$ such that $H\cdot b_2 = \phi(b_1,C^k)$ is the 
trivial subgroup.  

From the analysis above, we see that 
\[\models \forall x,x_1,x_2,{y}_1,  {y}_2 \, 
(\phi(x, x_1,x_2,{y}_1,  {y}_2)\rightarrow \displaystyle\bigwedge_{i,j} g_{ij}(y_{1j}) 
= y_{2i}).\]
Which implies 
\[\models \forall x,x_1,x_2,{y}_1\,\exists ! y_2\,(\phi(x,x_1,x_2,{y}_1,  {y}_2)).\]

With this observation, Corollary~\ref{cor:cnct-fn-1} yields that there is a 
$\emptyset$-definable $f:D\to C^k$ given by $\phi(x,x_1,x_2,{y}_1,  {y}_2)$ 
and $D = M^{t}\times C^n$ where $t = n+k+|x|$.  Since $D$ and $C^k$ are 
pre-smooth, and $f$ has a closed graph, by the Closed Graph Theorem 
(see~\cite[Lemma 5.5]{HZ96}) $f$ is continuous. Hence a morphism. 

For further reduction, one can take $D^0 := \{s_bm_1m_2\}\times \pr^{-1}(m_1)$, 
which is a finite, and hence a closed subset of $M^{t}\times C^n$. And directly 
show that $f\times\mathrm{id}: D^0\times C^m\to C^n\times C^m$ is 
continuous for all $n$. Hence, again a morphism.

Furthermore, as it is 
evident from the above analysis, the parameters $s_b,m_1,m_2$ are not 
essential. As we have shown $Mb_1$-locus of $b_2$ is the same as the 
locus of $b_2$ over $b_1$. But the $b_1$-locus of $b_2$ is given by 
$\displaystyle\bigwedge_{i,j} g_{ij}(y_{1j}) = y_{2i}$. Which defines a function 
\begin{align*}
C^n &\to C^k\\
b_1 &\mapsto b_2
\end{align*}
where $b_{2i}$ are obtained from the corresponding $b_{1j}$ and $g_{ij}$.

\begin{remark} 
As the parameter $s_b$ is not needed in this example, it seems 
that the function $f^\dagger$ will not play a role. Indeed, for this structure the 
analysis of the relevant relations of the type $g_{ij}(y_{1j}) = y_{2i}$ is enough to 
show that the specialisation extends from the base without using $f^\dagger$. 
For a proof see~\cite[Thm. 4.4.7]{Efem19}.
\end{remark} 

Although the function $f^\dagger$ does not play a role in this 
structure, we will still show that there is such a function so that 
this structure satisfies the Continuous Connections assumption. 
Fix an element $s_b\in M$. Although we established that it is much 
simpler, we may still view the function $f$ as taking any $(s_b,m_1,
m_2,\beta_1)\in D_0$ to a unique $\beta_2\in\pr^{-1}(m_2)$ for 
the sake of argument. Define $f^\dagger:D^\dagger\to M$ as 
$f^\dagger (m_1,m_2,\beta_1,\beta_2) = s_b$ where $D^\dagger 
= \{m_1m_2\}\times \pr^{-1}(m_1m_2)$. As $f^\dagger$ is a 
constant function it is clearly a morphism.  

Let $(s_b,m_1,m_2,\beta_1)\in D_0$, then $f(s_b,m_1,m_2,\beta_1) 
= \beta_2$ where $\beta_1\in\pr^{-1}(m_1)$ and $\beta_2\in
\pr^{-1}(m_2)$. Hence $(m_1,m_2,\beta_1,\beta_2)\in D^\dagger$, 
and $f^\dagger(m_1,m_2,\beta_1,\beta_2) = s_b$.



\end{document}